\documentclass{jfm} 

\usepackage{graphicx}
\usepackage{natbib}

\ifCUPmtlplainloaded \else
  \checkfont{eurm10}
  \iffontfound
    \IfFileExists{upmath.sty}
      {\typeout{^^JFound AMS Euler Roman fonts on the system,
                   using the 'upmath' package.^^J}%
       \usepackage{upmath}}
      {\typeout{^^JFound AMS Euler Roman fonts on the system, but you
                   dont seem to have the}%
       \typeout{'upmath' package installed. JFM.cls can take advantage
                 of these fonts,^^Jif you use 'upmath' package.^^J}%
      }
  \else
  \fi
\fi

\ifCUPmtlplainloaded \else
  \checkfont{msam10}
  \iffontfound
    \IfFileExists{amssymb.sty}
      {\typeout{^^JFound AMS Symbol fonts on the system, using the
                'amssymb' package.^^J}%
       \usepackage{amssymb}%
         \let\leq=\leqslant
         \let\geq=\geqslant
      }{}
  \fi
\fi

\ifCUPmtlplainloaded \else
  \IfFileExists{amsbsy.sty}
    {\typeout{^^JFound the 'amsbsy' package on the system, using it.^^J}%
     \usepackage{amsbsy}}
    {\providecommand\boldsymbol[1]{\mbox{\boldmath $##1$}}}
\fi

\newsavebox{\astrutbox}
\sbox{\astrutbox}{\rule[-5pt]{0pt}{20pt}}

\newcommand\p{\ensuremath{\partial}}

\usepackage{fixmath,amsmath,amsfonts}             
\usepackage{booktabs} 

\usepackage{hyperref}
\usepackage{amssymb}
\usepackage[labelfont=bf]{caption}
\usepackage{subcaption}
\usepackage{color}
\usepackage{float}

\newtheorem{remark}{Remark}

\newtheorem{problem}{Problem}

\def\div{{\rm div}}
\def\ep{\varepsilon}
\def\Oe{\Omega^\ep}
\def\ve{\mathbf{v}^\ep}
\def\o{\omega}
\def\O{\Omega}
\def\nab{\nabla}
\def\vphi{\boldsymbol\varphi}
\def\vphih{\boldsymbol\varphi_h}
\def\veh{\mathbf{v}^\ep_h}

\newcommand{\clos}[1]{\smash[t]{\overline{#1}}}

\newcommand{\cbl}{\ensuremath{C_1^{{bl}}}}
\newcommand{\cblh}{\ensuremath{C_{1,h}^{{bl}}}}
\newcommand{\cblw}{\ensuremath{C_\omega^{{bl}}}}
\newcommand{\cblwh}{\ensuremath{C_{\omega,h}^{{bl}}}}

\newcommand{\grid}{\ensuremath{{\cal T}_{h}}}

\newcommand{\Omegaf} {\ensuremath{{\Omega_{\operatorname{f}}}}}

\newcommand{\Omegap} {\ensuremath{{\Omega_{\operatorname{p}}}}}

\newcommand{\Gammaper} {\ensuremath{{\Gamma_{\operatorname{per}}}}}

\providecommand{\norm}[1]{\lVert#1\rVert}
\newcommand{\abs}[1]{\lvert#1\rvert}

\newcommand{\Set}[1]{\ensuremath{\left\{#1\right\}}}

\newcommand{\betabl}{\beta^{bl}}

\newcommand{\wbl}{\omega^{bl}}

\newcommand{\eff}{\mathit{eff}}
\newcommand{\IF}{\Sigma}
\newcommand{\per}{\mathit{per}}

\newcommand{\Zbl}{\ensuremath{Z^{bl}}}
\newcommand{\Zkl}{\ensuremath{Z^{k}_{l}}}
\newcommand{\cell}{\ensuremath{T}}
\newcommand{\transf}{\ensuremath{\Pi}}

\title[Pressure jump interface law]{Pressure jump interface law for the Stokes-Darcy coupling: Confirmation by direct numerical simulations}

\author[T. Carraro and C. Goll and A. Marciniak-Czochra and A. Mikeli\'c]%
{C\ls a\ls r\ls r\ls a\ls r\ls o\ns T.$^1$, \ns
G\ls o\ls l\ls l\ns  C.$^1$, \ns
M\ls a\ls r\ls c\ls i\ls n\ls i\ls a\ls k-C\ls z\ls o\ls c\ls h\ls r\ls a\ns A.$^{1,2}$\ns
\and M\ls i\ls k\ls e\ls l\ls i\ls \'c\ns A.$^3$
\thanks{Email address for correspondence: {\tt andro.mikelic@univ-lyon1.fr}.
}}

\affiliation{$^1$Institute for Applied Mathematics, Heidelberg University, 69120
Heidelberg, Germany\\[\affilskip]
$^2$Bioquant, Heidelberg University, 69120 Heidelberg, Germany\\[\affilskip]
$^3$Universit\'e de Lyon, CNRS UMR 5208,\\
  Universit\'e
Lyon 1, Institut Camille Jordan, \\   43, blvd. du 11 novembre 1918,
 69622 Villeurbanne Cedex, France
}

\date{Version from \today}
\begin{document}

\maketitle

\begin{abstract}
It is generally accepted that the effective velocity of a viscous flow over
a porous bed satisfies the Beavers-Joseph slip law. To the contrary,
interface law for the effective stress has been a subject of controversy.
Recently, a pressure jump interface law has been rigorously derived by
Marciniak-Czochra and Mikeli\'c.  In this paper, we provide a confirmation of the analytical result using direct numerical simulation of the flow at the microscopic level.
\end{abstract}

\begin{keywords}
interface pressure jump law, Beavers and Joseph law, Stokes-Darcy coupling, adaptive finite element, microscopic flow simulation
\end{keywords}
\section{Introduction}
\label{intro}
The focus of our paper is in mathematical description of slow incompressible viscous flow over a porous bed.  Physically, there is no interface between the unconstrained flow domain and the pores part of the porous bed.  For computational purposes we upscale the Navier-Stokes equation in the porous bed and replace it by
Darcy law. As opposed to the pore scale equations, which are defined in the pore structure, Darcy law is valid
at every point of the porous bed. Mathematically, the two systems of equations are very different: In the unconstrained flow domain description involves the Stokes
system. To describe the porous bed Darcy system is used.  Orders of the differential operators acting on velocity and
on pressure are different in these two cases.

After upscaling, an interface appears between the two domains and  relevant boundary conditions at the interface have
to be defined. We recall that the interface is not physical. It results
from the upscaling procedure and is needed to solve the effective equations. We can say
it is a computational interface.

Continuity of pressure and normal velocity are generally accepted interface conditions. However, they deserve some
comments. Darcy law for slow viscous flow through rigid porous medium was derived using a two-scale technique in ref.
 \cite{EneS:1975},  using rigorous homogenization for periodic porous media in ref. \cite{Tartar:1980}, and for random
porous media in ref. \cite{BeliaK:1995}. Derivation using volume averaging can be found in ref. \cite{Whitaker:1986}. All these derivations
require spatial homogeneity of the porous medium. They are in general not valid if the homogeneity is broken by, for example,
 presence of an interface with another medium.

Upscaled velocity in a porous medium remains incompressible. Global incompressibility implies continuity of
normal velocities at the interfaces. For the pressure field the situation is different. Physically, we have continuity
of the contact forces at any fluid interface. Calculating contact forces involves stresses. Upscaling of the
viscous stress in a rigid porous medium gives the effective pressure as the leading order stress approximation. Close to
interfaces and closed boundaries, the influence of the viscous part of the stress tensor might be much more
important, as we will see in the discussion which follows.

For a flow over a porous bed, another condition is added. It is the slip law published in ref. \cite{BeaveJ:1967}. In two dimensions (2D),
for a flat interface free fluid/porous bed, placed on $\{ x_2 =0 \}$, it reads
\begin{equation}\label{AUX}
   v_1 -v_D = \frac{\sqrt{k}}{\alpha} \frac{\partial v_{1}}{ \partial x_2} ,
\end{equation}
where $x_2$ is the direction orthogonal to the interface, $v_1$ is the free fluid tangential velocity at $\{ x_2 =0 \}$
and $v_D$ is the tangential component of Darcy velocity at the interface.
$k$ is  (the scalar) permeability and $\alpha$
characterizes the geometrical structure of the pores, close to the interface. The slip law (\ref{AUX}) was introduced on the
experimental basis by \cite{BeaveJ:1967}, for the case of a 2D channel Poiseuille flow over a porous bed.
Reference \cite{Saffman:1971} used Slattery's ad hoc form of Darcy law to give an analytic justification and remarked that
$v_D =O(k)$ and it can be dropped from (\ref{AUX}). Similar approach was undertaken in ref. \cite{Dagan:1981}.

More detailed analysis of the interface behavior is found in the seminal paper \cite{EneS:1975}. The authors
study interface conditions when the size of the velocity and the pressure are of the same order on both sides
of the interface. A dimensional analysis, and then boundary layer matching, were undertaken to conclude the pressure
continuity up to the order of pore size. Their considerations did not confirm the law by Beavers and Joseph. In ref \cite{LevyS:1975} further dimensional analysis arguments are provided in favor of pressure continuity.

Stationary incompressible viscous flow, described by the Stokes equations, adjunct to a 2D periodic porous medium,
was considered by ref. \cite{JaegerM:1996}, where also appropriate boundary layers were constructed. In the setting of
the experiment by ref. \cite{BeaveJ:1967}, reference \cite{JaegerM:2000} derived rigorously the slip law (\ref{AUX}) and calculated
the coefficient $\alpha$ and the permeability tensor $k$. For a flow over a porous bed, the effective velocity and
the pressure in the unconstrained fluid domain are calculated by imposing a non-penetration condition and the
Beavers-Joseph condition (\ref{AUX}).  Accuracy of the Darcy approximation in a porous medium is expressed using the dimensionless small parameter $\ep$, being the ratio of the typical pore size and the domain size.   In  ref. \cite{JaegerM:2000} the difference between the physical and effective (upscaled)
velocities and pressures was estimated in the energy norm, which involves squared volume integrals of velocity difference, 
velocity gradient difference and pressure difference
squared, by a constant times appropriate power of
the homogenization parameter $\varepsilon$. This error estimate,
providing a rigorous justification of the Beavers and Joseph law, is based on the construction of the corresponding boundary
layer around the interface. The slip coefficient is proportional to the viscous energy of the boundary layer. Numerical calculations of the slip coefficient, using the solution of the boundary layer equations for both isotropic and anisotropic
porous media, are found in ref. \cite{JaegerMN:2001}.

Prior to the analytical justification of law by Beavers and Joseph, several authors have attempted a numerical
justification of the slip law, by directly solving the Navier-Stokes equations in the unconstrained fluid and in
the porous bed.

References \cite{LarsoH:1986, LarsoH:1987} include simulations of the Stokes flow through an array of infinite, periodically distributed, cylinders.
The array is placed in the lower half-space and leads to a 2D problem. The authors placed the interface at the plane
 passing through the first row of cylinders. Due to this placement, they experienced difficulties in justifying the Beavers and
Joseph slip law and concluded that it is sensitive to the choice of the interface position.

In the study \cite{SahraK:1992}, the non-stationary Navier-Stokes equations were solved in both domains.
Simulations provided a justification of the Beavers and Joseph interface law. Dependence of
the slip condition on the position of the interface was studied as well and Saffman's predictions were confirmed. For
large Reynolds numbers, symmetry breaking affecting the value of the slip coefficient was observed. These
results are also discussed in the book \cite{Kaviany:1995}.

More recently ref. \cite{LiuP:2011} studied the finite Reynolds number 3D flow in a channel
with porous lateral walls by direct simulation, and concluded that the slip coefficient depends on the Reynolds number.

There exist important works in literature in which the coupling between free fluid flow and a porous medium is modeled by Brinkman
equation (the single domain approach). As remarked in ref. \cite{Nield:2009},
the model is semi-empirical and it gives satisfactory results only when two different viscosities are used in different places of the model. We will limit our attention to Darcy law and to the two-domain modeling.
Let us remark that since Brinkman equation is of the same type as the Navier-Stokes system, there is more
liberty in setting interface conditions. 
We mention in this context the model of ref. \cite{OchoaW:1995a, OchoaW:1995b}, where
the velocity is continuous and a stress jump is imposed.

Recapitulating, we have two interface conditions:
\begin{itemize}
  \item Continuity of the normal velocities.
  \item The slip law (\ref{AUX}) of Beavers and Joseph. 
\end{itemize}
Condition for the
normal stress at the interface has to be defined.
Generally accepted condition seems to be continuity of pressure. Since pressure
in the porous bed is an average, this law does not follow from continuity
of contact forces (the third Newton law).
Pressure continuity was strongly advocated in ref. \cite{EneS:1975}. The authors there propose a boundary layer matching
argument, but their argument is also used to contradict the slip law (\ref{AUX}) of Beavers and Joseph. On the other
 hand, direct numerical solution of the Navier-Stokes equations in ref. \cite{SahraK:1992} indicate
 a pressure slip for large Reynolds numbers.

Analytical answer is provided in the recent article \cite{MarciM:2012}, where the pressure interface
 condition for the Stokes flow over a porous bed was studied. Using homogenization coupled with the boundary layer
analysis around the interface, the authors obtained a rigorous pressure jump formula at the interface
\begin{gather}
 -[p] = p_D (x_1 , 0- ) - p_f (x_1 , 0+ ) =  {\cblw{}}  \frac{\p v_1}{ \p x_2 } (x_1 , 0 )  .  \label{Presspm2}
\end{gather}
$p_f$ is the effective unconstrained fluid pressure, $p_D$ the porous
medium pressure and $\cblw{}$  is a constant depending only on the porous bed geometry.
It corresponds to the pressure drop between $+\infty$ and $-\infty$, in the boundary layer problem used to calculate
the slip constant.

The law \eqref{Presspm2} was already proposed in ref. \cite*{JaegerMN:2001}.  More precisely,
the boundary layer equations, for the slip constant were numerically solved. For isotropic porous media it was
proved that {$\cblw{} =0$}, by a symmetry argument. Numerical simulations involving anisotropic periodic porous layers result in
{the order of} $\cblw{}= O(1)$ with respect to
the characteristic size of the pore. This leads to the conclusion that pressure is not continuous in general.

For the Navier-Stokes flow over a porous bed anisotropy is due to the inertia term and ref. \cite{SahraK:1992}
also observed numerically, for increasing Reynolds number, the pressure slip.


The goal of the present paper is to confirm the pressure jump law \eqref{Presspm2} by a direct simulation of the Stokes flow
in a porous bed and  unconstrained domain, corresponding to the
Beavers and Joseph experiment.
This provides,
besides the analytical proof from ref. \cite{MarciM:2012}, another independent confirmation by direct
simulation. We will also justify the law (\ref{AUX}) by Beavers and Joseph. Also we note that dependence
of the slip coefficient on the change of position of the interface, has been studied analytically in detail in ref. \cite{MarciM:2012}, and affine dependence was found.

The outline of the paper is as follows: In Section~\ref{problem_setting} we present an asymptotic expansion which yields both the slip law
of Beavers and Joseph and the pressure jump law (\ref{Presspm2}). Asymptotic expansion was rigorously justified in ref.
\cite{MarciM:2012}. 
In Section~\ref{numerical_setting} we present the adaptive finite element
approach to solve the problems needed for the numerical comparison between
microscopic and upscaled problems.
Section \ref{numerical_confirmation} shows results of numerical simulations for two
flow cases confirming the laws (\ref{AUX}) and (\ref{Presspm2}).

\section{Problem setting}
\label{problem_setting}
 We assume a slow incompressible viscous flow through an unconfined region $\Omega_f=(0,L)\times(0,1)$, where $L$ is a positive number denoting the length of the domain, and the fluid part of a porous medium  $\Omega_p=(0,L)\times(-1,0)$. The domain is assumed to be sufficiently large and, after subtracting the boundary conditions, we are allowed to consider periodicity in the longitudinal direction. We start by setting the geometry and writing the dimensionless flow equations.
We assume the no-slip condition at the boundaries of the pores (i.e., a rigid porous medium).
\subsection{Microscopic equations}\label{Microgeoeq}
\begin{figure}
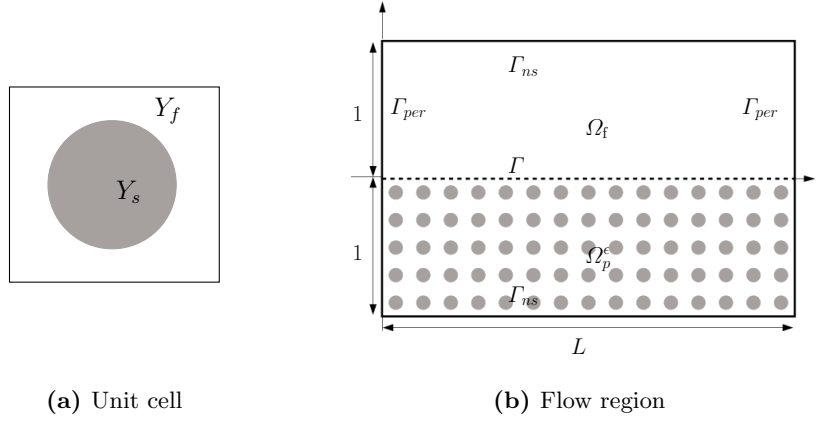

\centering
\begin{subfigure}[b]{0.3\textwidth}
\centering
\raisebox{1.2cm}{\resizebox{0.7\textwidth}{!}{\input{UnitCell}}}
\caption{Unit cell}
\label{subfig.unitcell}
\end{subfigure}
\begin{subfigure}[b]{0.6\textwidth}
\centering
\resizebox{0.8\textwidth}{!}{\input{Geometry_per}}
\caption{Flow region}
\label{TwoDomains}
\end{subfigure}
\caption{The model geometry}
\end{figure}

The geometry of the problem is given in figure~\ref{TwoDomains} and, more precisely, the periodic structure is defined as follows.
The porous part has a periodic structure and it corresponds to a repetition
of so-called  cells of the characteristic size $\ep$. Each cell
is made from the unit cell $Y=(0,1)^2$, rescaled by $\ep$. The unit cell contains a
pore part $Y_f$ and a solid part $Y_s$ with $Y_s \subsetneq Y$ (see
figure~\ref{subfig.unitcell}). The union of all pores gives the fluid part
$\Omega_p^\ep$  of porous domain $\Omega_p$.   $\Gamma=(0,1)\times \{0\}$
describes an interface between the unconfined domain and the porous domain.
The flow takes place in $\Omega^\ep=\Omega_p^\ep \cup \Omega_f \cup
\Gamma$
 and it is described by the following non-dimensional steady Stokes system
 in $\Omega^\varepsilon $:
 \begin{gather} -
 \upDelta \mathbf{v}^\ep + \nabla p^{\ep} = \mathbf{f} \qquad \hbox{ in } \quad \Oe
\label{1.3} \\ \div \, \ve = 0 \qquad \hbox{ in } \quad \Oe , \qquad \int_{\O_f} p^\ep \ dx =0,
\label{1.4} \\ \ve  =0 \; \hbox{on } \ \p \Oe \setminus \bigg(  \{ x_1 = 0 \} \cup \{ x_1 = L \} \bigg)
, \qquad \{ \ve , p^\ep \} \ \hbox{ is }
L-\hbox{periodic in } \; x_1 . \label{1.5} \end{gather}
  Here the non-dimensional  $\mathbf{f}$ stands for the effects of external forces or an
injection at the boundary or a given pressure drop, and it corresponds to the physical forcing term
multiplied by the ratio between Reynolds' number and Froude's number squared.  Specifically, if the force $\mathbf f$ is
non-constant, it corresponds to a non-constant pressure drop or to a
non-parabolic injection profile.
$\ve $ denotes the non-dimensional velocity and  $p^\ep$ is the non-dimensional
pressure.

\subsection{Two-scale expansion}

 The idea behind the two-scale expansion is the following: without forcing
 infiltration, in the interior of the porous medium, the permeability is
 $k= O(\ep^2)$ and Darcy's velocity is small. Consequently, the flow is tangent to the interface $\Gamma$. The leading order approximation of system (\ref{1.3})-(\ref{1.5}) is the Stokes flow in $\Omega_f$, with no-slip condition on $\Gamma$. It is modeled by the system,
\begin{gather}
- \upDelta \mathbf{v}^0  + \nabla p^0 = \mathbf{f}
\qquad \hbox{ in } \O_f ,\label{4.37}\\
\div \  \mathbf{v}^0 = 0 \qquad \hbox{ in } \O_f , \qquad \int_{\O_f} p^0 \ dx =0, \label{4.38}\\
\mathbf{v}^0 = 0  \quad  \hbox{ on } \p   \O_f  \setminus \bigg(  \{ x_1 = 0 \} \cup \{ x_1 =L \} \bigg)
, \quad
\{ \mathbf{v}^0 , p^0 \}    \quad  \hbox{ is } \;
L\hbox{-periodic in } \; x_1. \label{4.40}
\end{gather}
Following the two-scale expansions from \cite{EneS:1975},  the behavior of  the velocity and pressure fields in $\O_p$, far from the outer boundaries, is expected to be given by the following system of equations,
\begin{gather}
\mathbf{v}^\ep (x) =  \ep^2 \sum_{j=1}^2 \mathbf{w}^j (\frac{x}{\ep}) ( f_j (x) -
\frac{\partial p_D (x)}{ \partial x_j } ) + O(\ep^3 )\quad x\in \O_p , \label{1.62} \\
  p^\ep (x) = p_D (x) +  \ep \sum_{j=1}^2 \pi^j (\frac{x}{\ep}) ( f_j (x) -
 \frac{\partial p_D (x)}{ \partial x_j } ) + O(\ep^2) \quad x\in \O_p ,\label{1.63} \\
  \mathbf{v}^D (x) = K
(\mathbf{f} (x) - \nabla_x p_D (x) )  \, (\hbox{\bf Darcy's law}) ,
\label{1.64} \\
\mbox{div } \mathbf{v}^D  =0 \quad \mbox{ on } \; \O_p, \label{1.64A}\end{gather}
where $\{ \mathbf{w}^j , \pi^j \}$ are calculated using cell problems and $K$ consists of the volume averages of $\mathbf{w}^j$, $j=1,2$, see Subsection~\ref{subsec.cell problem}. Note that the dimensionless permeability is $\ep^2 K$ and it is a symmetric positive definite matrix. See e.g. \cite{JaegerM:2009} or \cite{Allaire:1997} for more details.

System (\ref{1.64})-(\ref{1.64A}) describes the effective pressure $p_D$. However, we do not know the boundary condition for $p_D$ (respectively $\mathbf{v}^D$) on the interface $\Gamma$, and $\{\mathbf{v}^D,  p_D\} $ are not determined.

We need \textbf{interface conditions} coupling problem (\ref{4.37})-(\ref{4.40}) with (\ref{1.64})-(\ref{1.64A}).

Natural approach to find interface conditions is by using matched asymptotic expansions (MMAE). The method is used in a number of situations arising in mechanics. 
  For a detailed presentation of the MMAE method we refer to the book  \cite{Zeytounian:2002}  and to references therein.

In the language of the MMAE, expansions in $\O_f$  and $\O_p$  are called the outer expansions. The boundary and/or interface behavior is captured by an inner expansion.  In the inner expansion  the independent variable is stretched out in order to describe the behavior in the neighborhood of the boundary and/or interface.

The MMAE approach  matches the two expansions. In the singular perturbation
problems involving boundaries, only the function values at the boundary are
matched and the approach works well. When interfaces are involved, it is
needed to match additionally the values of the normal derivatives.
This difficulty is not easy to circumvent because imposing matching of the values of the function and its normal derivative leads to an ill posed problem for the second order equation.
This difficulty can not been easily avoided, since the simultaneous imposition of matching conditions for the values and
for the normal
derivative of the function leads to an ill posed problem for the second
order equation.


Here, Darcy's velocity in $\O_p$ is of order $O(\ep^2)$. Therefore at the lowest order MMAE confirms the boundary condition (\ref{4.40}) on $\Gamma$, i.e. $\mathbf{v}^0 =0$.  
 Additional physical matching conditions would be the continuity of the
 contact forces, which is not
 assured by MMAE. Therefore, it is not clear if we are allowed to match
 the values of the pressure on $\Gamma$.

The absence of a matching condition in the velocity gradient and in the pressure leads to a jump of the contact force.




 The difficulty was solved  using a boundary layer correction in \cite{JaegerM:2000},  \cite*{JaegerMN:2001} and \cite{MarciM:2012}. The applied strategy is the following:
 We handle the pressure jump  by adjusting the porous medium pressure and
 the shear stress jump using a particular {properly
 derived} boundary layer. In fact, the shear stress jump influences the pressure values as well.

\begin{figure}
\center
\resizebox{0.2\textwidth}{!}{\input{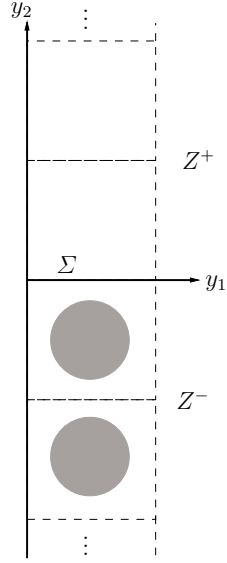}}
\caption{Domain of the Navier boundary layer problem.}
\label{fig.boundary_layer_domain}
\end{figure}

At the interface  $\Gamma$  we have the shear stress jump equal to $\displaystyle - {\partial
v_1^{0}}/{\partial x_2}|_{\Gamma}$. The natural stretching
variable is given by the geometry and it coincides with the fast variable   $\displaystyle
y={x}/{\ep}$. The correction $\{ \mathbf{w} , p_w \}$ to the zero order approximation satisfies again the Stokes system
\begin{gather}  \bigl[ \mathbf{w} \bigr] (y_1 , 0)= \mathbf{w} (y_1 , 0+) - \mathbf{w} (y_1 , 0-) =0, \quad [p_w ] (y_1 , 0) =0\notag \\
  \mbox{ and } \quad \bigr[  \frac{\partial w_1}{\partial y_2}
\bigl] (y_1 , 0) =   \frac{\partial v_1^{0}}{\partial x_2} (x_1 ,0)
|_{\Gamma} \quad \hbox{ on }  \quad
{\frac{\Gamma}{\ep}} .\label{BJ5.4)}
\end{gather}
Using periodicity of the geometry and independence of $\displaystyle
\frac{\partial v_1^{0}}{\partial x_2} |_{\Gamma}$ of the fast variable $y$, we obtain
\begin{equation}\label{BJ5.6}
    \mathbf{w} (y) =  \frac{\partial
v_1^{0}}{\partial x_2} |_{\Gamma} \boldsymbol\beta^{bl} (y) \quad \mbox{ and }
\quad p_w (y) =  \frac{\partial v_1^{0}}{\partial x_2} |_{\Gamma}
\omega^{bl} (y),
\end{equation}
where $\{ \boldsymbol\beta^{bl} , \omega^{bl} \}$ is calculated  in a semi-porous column $Z_{BL} = Z^+ \cup \Sigma \cup Z^- $,  with $\Sigma=(0,1)\times \{ 0\} $,  $Z^+ = (0,1) \times (0, +\infty )$,
and the semi-infinite porous slab is $Z^- =
\cup_{k=1}^\infty ( Y_f -\{ 0,k \} )$. See figure~\ref{fig.boundary_layer_domain} for more details.

If $D_y$ denotes the symmetrized gradient, then $\{ \boldsymbol\beta^{bl} , \omega^{bl} \}$ is given by
\begin{gather}
-\upDelta _y \boldsymbol\beta^{bl} +\nabla_y \omega ^{bl} =0\qquad \hbox{ in
} Z^+ \cup Z^-, \label{BJ4.2}\\ \div_y  \boldsymbol\beta^{bl} =0\qquad \hbox{ in }
Z^+ \cup Z^-, \label{4.3} \\ \bigl[ \boldsymbol\beta^{bl} \bigr]_\Sigma (\cdot , 0)= 0
 \quad \mbox{ and } \quad \bigr[ \{ 2 D_y (\boldsymbol\beta^{bl} )
-\omega^{bl} I \}  \mathbf{e}^2 \bigl]_\Sigma (\cdot , 0) = \mathbf{e}^1 \quad \hbox{ on }
\Sigma \label{4.5)}, \\\boldsymbol \beta^{bl} =0 \quad \hbox{ on }
\displaystyle{\bigcup_{k=1}^{\infty} ( \p Y_s} -\{ 0,k \} ), \qquad \{
\boldsymbol\beta^{bl} , \omega^{bl} \} \quad \hbox{ is } 1\hbox{-periodic in }
y_1. \label{4.6}
\end{gather}
The problem (\ref{BJ4.2})-(\ref{4.6}) was studied in \cite{JaegerM:1996} and it was proved that
\begin{itemize}
  \item Gradients of $\{ \boldsymbol\beta^{bl} , \omega^{bl} \}$ stabilize exponentially fast to $0$.
  \item $ \boldsymbol\beta^{bl} $  stabilizes exponentially fast to $C_1^{bl} \mathbf{e}^1$, when $y_2 \to +\infty$ and to zero when $y_2 \to -\infty$. $C_1^{bl} $ is strictly negative.
  \item $\omega^{bl} $ stabilizes exponentially fast to $C_\o^{bl} \mathbf{e}^1$, when $y_2 \to +\infty$ and to $C^{bl}_{0}$ when $y_2 \to -\infty$. Since we have liberty in adding a constant to the pressure, we choose $C^{bl}_{0}=0$.
\end{itemize}
In fact, it is the absence of stabilization of the boundary layer velocity $\boldsymbol\beta^{bl}$ to zero which yields a slip.
In addition, $\mathbf{w}$ can't be a correction
because of the  stabilization of $\boldsymbol\beta^{bl  }$
towards a nonzero constant velocity $  C^{bl}_1 \mathbf{e}^1 $. It creates a counterflow
at the upper boundary of $\O_f$,  given by the following
 Stokes system  in $\Omega_f$:
\begin{gather}
- \upDelta \mathbf{z}^\sigma  + \nabla p^\sigma = 0
\qquad \hbox{ in } \O_f ,\label{4.37Couette}\\
\div \  \mathbf{z}^\sigma = 0 \qquad \hbox{ in } \O_f ,\label{4.38Couette}\\
\mathbf{z}^\sigma = 0  \quad  \hbox{ on }  \{ x_2 = 1 \} \quad \mbox{and} \; \mathbf{z}^\sigma =\frac{\p   v^0_1 }{ \p x_2 }
|_\Gamma  \mathbf{e}^1 \quad  \hbox{ on }  \{ x_2 = 0 \}
,\label{4.39Couette}\\ \{ \mathbf{z}^\sigma , p^\sigma \}    \qquad  \hbox{ is } \;
1\hbox{-periodic in } \; x_1 .\label{4.40Couette}
\end{gather}
Now we are in the situation to propose  the two-scale expansion for the velocity:
\begin{gather}
\ve = \underbrace{\mathbf{v}^0 -   \ep  C^{bl}_1  \mathbf{z}^\sigma
}_{\mbox{the outer expansion}}  \underbrace{-  \ep (\boldsymbol\beta^{ bl} (\frac{x}{ \ep}) -   C^{bl}_1 \mathbf{e}^1 )  \frac{\p   v^0_1 }{ \p x_2 }
|_\Gamma}_{\mbox{the inner expansion}} +\dots \quad \mbox{in } \; \O_f ,
 \label{4.66} \\
 \ve = \underbrace{O(\ep^2)}_{\mbox{the outer expansion}}  \underbrace{-
 \ep \boldsymbol\beta^{ bl } (\frac{x}{ \ep})  \frac{\p   v^0_1 }{ \p x_2 }
|_\Gamma}_{\mbox{the inner expansion}} +\dots \quad \mbox{in } \; \O_p .
 \label{4.66A}
\end{gather}
For the two-scale expansions (\ref{4.66})-(\ref{4.66A}) the values at the
interface $\Gamma$ are matched exactly and the shear stresses are matched
{with an approximation of order
$O(\epsilon)$}. At the flat interface $\Gamma$, with no slip condition for $\mathbf{v}^0$ and interface continuity of the boundary layer velocity, continuity of the normal component of the normal stress (i.e. of the normal contact force) reduces to the pressure continuity.

We need the two-scale expansion for the pressure.  Stabilization of the boundary layer pressure to $C^{bl}_\omega $,
when $y_2 \to +\infty$, influences strongly the pressure approximation. It reads
\begin{gather}
p^\ep  =\underbrace{p^0 {-}\> \ep  C^{bl}_1  p^\sigma}_{\mbox{the outer expansion}}  \underbrace{  - \bigl( {\omega}^{ bl} (\frac{x}{ \ep})  -
C^{bl}_\omega \bigr)   \frac{\p   v^0_1 }{ \p    x_2 } |_\Gamma}_{\mbox{the inner expansion}} +\dots \quad \mbox{in } \; \O_f , \label{BJ47} \\
 p^\ep  =\underbrace{p_D +  \ep \sum_{j=1}^2 \pi^j (\frac{x}{\ep}) ( f_j (x) -
 \frac{\partial p_D (x)}{ \partial x_j } ) }_{\mbox{the outer expansion}}  \underbrace{  - {\omega}^{ bl} (\frac{x}{ \ep})  \frac{\p   v^0_1 }{ \p    x_2 } |_\Gamma}_{\mbox{the inner expansion}} +\dots \quad \mbox{in } \; \O_p . \label{BJ47A}\notag
\end{gather}
Two-scale expansions (\ref{BJ47})-(\ref{BJ47A}) match at the interface $\Gamma$ at order $O(\ep)$ if and only if
\begin{equation}\label{BJ46}
     p^0 (x_1, +0) -  {p}_D (x_1, -0) =-  C^{bl}_\omega
  \frac{\p   v^0_1 }{ \p    x_2 } |_\Gamma   \quad \mbox{ for }
\quad x_1 \in (0,1).
\end{equation}
The conditions (\ref{BJ46}) allows to calculate Darcy's pressure $p_D$. It satisfies the equations (\ref{1.64})-(\ref{1.64A}), the condition (\ref{BJ46}), $v^D_2 =0$ on $\{ x_2 =-1 \}$ and $p_D$ is $1$-periodic in $x_1$.

{The results from \cite{MarciM:2012} (and from \cite{JaegerM:2000} in the case of Poiseuille's flow) yield in $\O_f \cup \O_p$
\begin{align*}
\ve -  \mathbf{v}^0  +\ep \left( \boldsymbol\beta^{ bl}(\frac{ x}{ \ep}) -   C^{bl}_1 \mathbf{e}^1 H(x_2)\right)  \frac{\p v^0_1
}{ \p   x_2 } |_\Gamma  +\ep   C^{bl}_1 \mathbf{z}^\sigma H(x_2 ) &=O (\ep^2), \\
p^\ep - p^0 H(x_2) - { p}_D  H(-x_2 ) +\bigl( {\omega}^{ bl , \ep} (x)
{-}\> H(x_2)C^{bl}_\omega \bigr)   \frac{\p   v^0_1 }{ \p    x_2 } |_\Gamma +
  \ep  C^{bl}_1  p^\sigma H(x_2)
   &=O (\ep ),    \\
  \nabla \ve -  \nabla \mathbf{v}^0  +\ep \nabla \bigg( ( \boldsymbol\beta^{ bl}(\frac{ x}{ \ep}) -   C^{bl}_1 \mathbf{e}^1 H(x_2))  \frac{\p v^0_1
}{ \p   x_2 } |_\Gamma  +  C^{bl}_1 \mathbf{z}^\sigma H(x_2 ) \bigg)&=O (\ep),
\end{align*}
where $H(t)$ is Heaviside's function, equal to one for $t>0$ and to zero for $t<0$. Furthermore, on $\Gamma$ it holds
\begin{align}
\ve  +\ep  \boldsymbol\beta^{ bl}(\frac{ x}{ \ep})   \frac{\p v^0_1
}{ \p   x_2 } |_\Gamma   =O (\ep^{3/2}).\label{equ.bjs epsilon32}
\end{align}}

\subsection{Interface condition}
\label{Interface condition}
The above results allow to find the effective interface conditions. Following \cite*{JaegerM:2000, JaegerM:2009, MarciM:2012}, we compare on the interface\/ $\Gamma$ the shear stress and the tangential velocity:
\[
\frac{\p v^\ep _1}{ \p x_2 } |_\Gamma = \frac{\p v^0 _1}{ \p x_2 }
|_\Gamma - \frac{\p \beta^{bl}_1 }{ \p y_2} |_{\Gamma , y=x/\ep}  +
O(\ep) \quad \hbox{ and\/ } \quad \frac{v^\ep _1}{ \ep } = -
\beta^{bl}_1 ( x_1 / \ep , 0) \frac{\p v^0 _1}{ \p x_2 } |_\Gamma +
O(\ep) .
\]
After averaging over\/ $\Gamma$ with respect to $y_1$, we obtain the 
Saffman version of the law by
Beavers and Joseph (\ref{AUX}), with
$$ \alpha = -\frac{\sqrt{k}}{\ep C_1^{bl}}, \quad k=O(\ep^2 ), \; C^{bl}_1 <0.$$
Therefore the effective flow in $\O_f$ is given by the following problem.
\begin{problem}[Effective flow]
Find a velocity field
$u^{\eff}$ and a pressure field $p^{\eff}$ such that
\begin{gather}
- \upDelta \mathbf{u}^{\eff}  + \nabla p^{\eff} = \mathbf{f} \qquad \hbox{ in } \O_f ,\label{4.91}\\
\div \ \mathbf{u}^{\eff} = 0 \qquad \hbox{ in } \O_f ,  \qquad \int_{\O_f} p^{\eff} \ dx =0,\label{4.92}\\ \mathbf{u}^{\eff} =
0  \qquad  \hbox{ on } (0,L) \times  \{1\}  
; \quad \mathbf{u}^{\eff} \; \hbox{ and }  \;  p^{\eff} \quad  \hbox{ are } \;
L\hbox{-periodic in} \quad x_1, \label{4.94}\\ u^{\eff}_2 = 0  \quad
\hbox{ and } \quad u^{\eff}_1 + \ep C^{bl}_1 \frac{\p u^{\eff}_1 }{ \p
x_2 } =0 \quad \hbox{ on } \quad \Gamma  . \label{4.95}
\end{gather}
\end{problem}
After \cite{MarciM:2012}, Theorem 2, we have

\begin{align}
\int_{\O_f} | \ve - \mathbf{u}^{\eff} |^2 \ dx   + \vert  M^\ep - M^{\eff}  \vert^2 &= O (
\ep^{3}), \label{4.100A} \\
\int_{\O_f} \{  |  p^\ep - p^{\eff}  |
+ |  \nabla (\ve - \mathbf{u}^{\eff} ) |\} \ dx &= O (\ep),
 \label{4.100}
\end{align}
where $\displaystyle M^{\eff} = \int_{\O_f} u_1^{\eff} \ dx $ is the mass flow.

The estimates (\ref{4.100A})-(\ref{4.100}), obtained analytically  in \cite{MarciM:2012}, will be verified by a direct numerical simulation in Section~\ref{numerical_confirmation}.

Next we recall that $p_D$ is given by
\begin{gather}
    \mbox{\rm div } \bigg( K  (\mathbf{f} (x) - \nabla p_D )\bigg) =0\; \mbox{ in\/ } \; \O_p,  \label{PresspmA} \\
 {p}_D  = p^{\eff} + C^{bl}_\omega  \frac{\p u^{\eff}_1}{ \p x_2 } (x_1 , 0+ )  \; \mbox{ on\/ } \; \Gamma; \quad
 K  (\mathbf{f} (x) - \nabla  { p}_D ) |_{\{ x_2 =-1 \} } \cdot \mathbf{e}^2 =0,  \label{Presspm2A}
\end{gather}
and after \cite{MarciM:2012}, Theorem 3, we have

\begin{gather}
    |\int_{\O_p} \{ \biggl( \ve + \ep \boldsymbol\beta^{bl} (\frac{x}{  \ep})  \frac{\p u^{\eff}_1}{ \p x_2 } (x_1 , 0 ) \biggr)  - K  (\mathbf{f} - \nabla  { p}_D ) \} \varphi \ dx |= o(1), \notag \\
     \mbox{ for every smooth } \; \varphi,  \quad \mbox{ as } \; \ep \to 0;   \label{ConcDarcy}\\
    \int_{\O_p} | p^\ep - { p}_D |^2 \ dx = o(1),   \quad \mbox{ as } \; \ep \to 0; \label{ConcPression} \\
    |\int_{\Gamma} ( p^\ep - p^{\eff} ) \varphi  \ dx_1 |  =O( \sqrt{\ep} ) \; \mbox{ for every smooth } \; \varphi,  \quad \mbox{ as } \; \ep \to 0.\label{EstPressBdry}
\end{gather}

The estimate (\ref{ConcPression}) will also be verified by a direct numerical simulation in Section~\ref{numerical_confirmation}.

\begin{remark}[Extension of velocity and pressure]
Fluid velocity $\mathbf v^\ep$ is extended zy zero to the solid part $\Omegap \setminus \clos{\Omega_p^\ep}$ of the porous medium $\Omegap$. If $Y^\ep_{f,i}$ is the $i$-th pore, then the pressure field $p^\ep$ is extended to the corresponding solid part $Y^\ep_{s,i}$ by setting
\begin{align}\label{equ.pressure extension}
 p^\ep(x) =
\begin{cases}
p^\ep,& x \in \O^\ep,\\
\frac{1}{\abs{ Y_{f,i}^\ep}}\int_{ Y_{f,i}^\ep} p^\ep,& x\in Y_{s,i}^\ep,
\end{cases}
\end{align}
where $\abs{ Y_{f,i}^\ep}$ denotes the volume of $ Y_{f,i}^\ep$.
The pressure extension \eqref{equ.pressure extension} is the extension of \cite{Lipton:A1989} and comes out from Tartar's construction, see \cite{Allaire:1997} for more details.
\end{remark}

\section{Finite Element Formulation}
\label{numerical_setting}
%
In this study the two conditions \eqref{4.95} and the first of \eqref{Presspm2A} are numerically confirmed by a direct simulation. 
To numerically verify all theoretical results we have to solve the following problems: the microscopic problem (\ref{1.3}--\ref{1.5}), the effective flow (\ref{4.91}--\ref{4.95}), the Darcy's law (\ref{PresspmA}--\ref{Presspm2A}), the boundary layer problem (\ref{BJ4.2}--\ref{4.6}) and some appropriate cell problems to compute the rescaled permeability $K$. 

All these problems have to be solved for two different kinds of inclusions.
The  solution of the microscopic problem has also to be computed  for
different boundary conditions, particularly a periodic configuration and a
flow with an injection boundary condition, see subsections \ref{sec:case_I} and \ref{sec:case_II}.
Particular attention has to be given to the calculation of the constants \cbl{} and \cblw{} used in the interface condition (see Section~\ref{Interface condition}), since we are going to show converge results with $\epsilon\to0$ in Section~\ref{numerical_confirmation} and the precision of several quantities reaches quickly the discretization error.

For this reason we adopt a goal oriented adaptive scheme for the grid refinement that allows reducing the computational costs to obtain a precise evaluation of a given functional, in particular to compute the two constants.

\subsection{Finite Element Formulation of the Microscopic Problem}
\label{Numset}
To numerically solve the problems we consider the finite element method and we give exemplary in this subsection the formulation of the microscopic problem \eqref{1.3}-\eqref{1.5}. For a more thorough introduction into the theory of finite elements, we refer to standard literature such as \cite{Ciarlet:2002} or \cite{BrennS:2002}.

The natural setting of the finite element approximation of the problem is its weak formulation, shown below.
We first introduce the spaces
\begin{align}
V(\Oe) &:= \Set{\mathbf{v} \in H^1(\Oe)^2\ \big|\ \mathbf{v}  =0 \text{ on }  \p \Omega \setminus \Gammaper, \mathbf{v} \hbox{ is } L\text{-periodic in } x_1},\\
L_0(\Oe) &:= \Set{ p \in L^2(\Oe)\ \big | \ \int_{\Oe} p \ dx = 0},
\end{align}
where $L^2(\Oe)$ is the space of square-integrable functions in $\Oe$, i.e. for $u \in L^2(\Oe)$ holds $\displaystyle \int_{\Oe} \abs{u(x)}^2\,dx < \infty$, and $H^1(\Oe)$ is the space of square-integrable functions, with first derivatives also square-integrable.

The weak formulation of problem \eqref{1.3}-\eqref{1.5} reads as follows:
\begin{problem}[Microscopic Problem in Weak Formulation]\label{prob.micro_weak}
For given $\mathbf{f}$ find a pair $(\ve,p^\ep) \in V(\Omega^\epsilon)\times L_0(\Oe)$, such that
\begin{align}
\int_{\Oe} \big(\nab \ve + (\nab \ve)^T\big) \cdot \nab \vphi \ dx + \int_{\Oe} p^\ep \nab \cdot \vphi \ dx &= \int_{\Oe} \mathbf{f} \cdot \vphi \ dx &&\forall \vphi \in V(\Oe),\\
\int_{\Oe} \nab \cdot \ve \ \psi \ dx &= 0 &&\forall \psi \in L_0(\Oe).
\end{align}
\end{problem}

We use finite elements to discretize this problem and consider a decomposition $\grid$ of the domain into so called cells $T$, whose union constitutes an approximation of the problem geometry, i.e. $\grid = \Set{\cell}$. 
We consider shape regular grids.
The cells are constructed via a set of polynomial transformations $\Set{\transf_\cell}_{\cell \in \grid}$ of a unit reference cell $\hat \cell$, see also Remark~\ref{rem.transformation}. The diameters $h_\cell$ of the cells define a mesh parameter $h$ by the piecewise constant function $h_{|\cell} = h_\cell$. 
On the grid we define for $s \in \mathbb N$
the finite dimensional space 
\begin{align}
\mathcal{S}_h^s(\Oe):=\Set{v_h\in C^0(\overline {\Oe})\ \big| \ {v_h}_{|\cell} \in Q^s(\cell), \cell \in \grid },
\end{align}
where $C^0(\overline \Oe)$ is the space of continuous functions on $\overline {\Oe}$. Let $P^s(\hat \cell)$ be space of polynomials of order lower or equal to $s$, then the space 
\begin{align}
Q^s(\cell)= \Set{p: T \to \mathbb{R}\ |\ p\left(\transf_\cell(\cdot)\right) \in  P^s(\hat \cell)}
\end{align}is the space of functions obtained by a transformation of bilinear ($s=1$), biquadratic ($s=2$) and in general higher order polynomials defined on the unit reference cell $\hat \cell$. 

For convergence results with respect to $h$ we consider a family of grids obtained by either uniform or local refinement of an initial regular grid.

\begin{remark}[Boundary Approximation]\label{rem.transformation}
Since the considered domains have curved boundaries we correspondingly use cells with curved boundaries (i.e. isoparametric finite elements) to get a better approximation. 
Considering the space $\mathcal{S}_h^2$ for the velocity we use biquadratic transformations of the unit cell $\hat \cell$.
\end{remark}
For the discretization of the Stokes system we use the Taylor-Hood element that uses the ansatz space $V_h(\Oe):=\left(\mathcal{S}^2_h(\Oe)\right)^2 \cap V(\Oe)$ for the velocity and $L_h(\Oe):=\mathcal{S}_h^1(\Oe)$ for the pressure.
This discretization is inf-sup stable (cf. \cite{BrezzF:1991}), so we do not need stabilization terms to solve the saddle point corresponding to the Stokes system, as for example in \cite*{JaegerMN:2001}.

The \textbf{finite element approximation} of the microscopic problem is obtained by replacing the (infinite dimensional) function spaces $V(\Oe)$ and $L_0(\Oe)$ by their discretized counterparts $V_h(\Oe)$ resp. $L_h(\Oe)$.
\begin{problem}[Finite Element Approximation of Microscopic Problem]\label{prob.micro_fe}
Find a pair $(\veh,p^\ep_h) \in V_h(\Omega^\epsilon)\times L_h(\Oe)$, such that
for all 
$(\vphih,\psi_h )  \in V_h(\Oe)  \times L_{0,h}(\Oe)$
\begin{align}
\int_{\Oe} \big(\nab \veh + (\nab \veh)^T \big) \cdot \nab \vphih \ dx + \int_{\Oe} p^\ep_h \nab \cdot \vphih \ dx &= \int_{\Oe} \mathbf{f} \cdot \vphih \ dx ,\\
\int_{\Oe} \nab \cdot \veh \ \psi_h \ dx &= 0 
\end{align}
and $\int_{\Oe} p^\ep_h\ dx = 0$.
\end{problem}

As shown in Section~\ref{problem_setting}, see also \cite{JaegerM:1996,JaegerM:2000,JaegerM:2009},  the Navier boundary layer problem (\ref{BJ4.2})-(\ref{4.6}) has to be solved to determine the constants \cbl{} and \cblw{} in the interface law \eqref{4.95} and the first of \eqref{Presspm2A}.
Next subsection is thus dedicated to the numerical determination of these constants and we will show in sections \ref{sec:case_I} and  \ref{sec:case_II} by direct numerical solving of the microscopic problem that the constant $\cblw$ is related to the pressure difference between the free fluid and the porous part.

\begin{figure}
\centering
 \begin{subfigure}[htb]{0.3\textwidth}
\centering 
\includegraphics[trim=54mm 80mm 21mm 83mm, clip, width =0.8\textwidth]{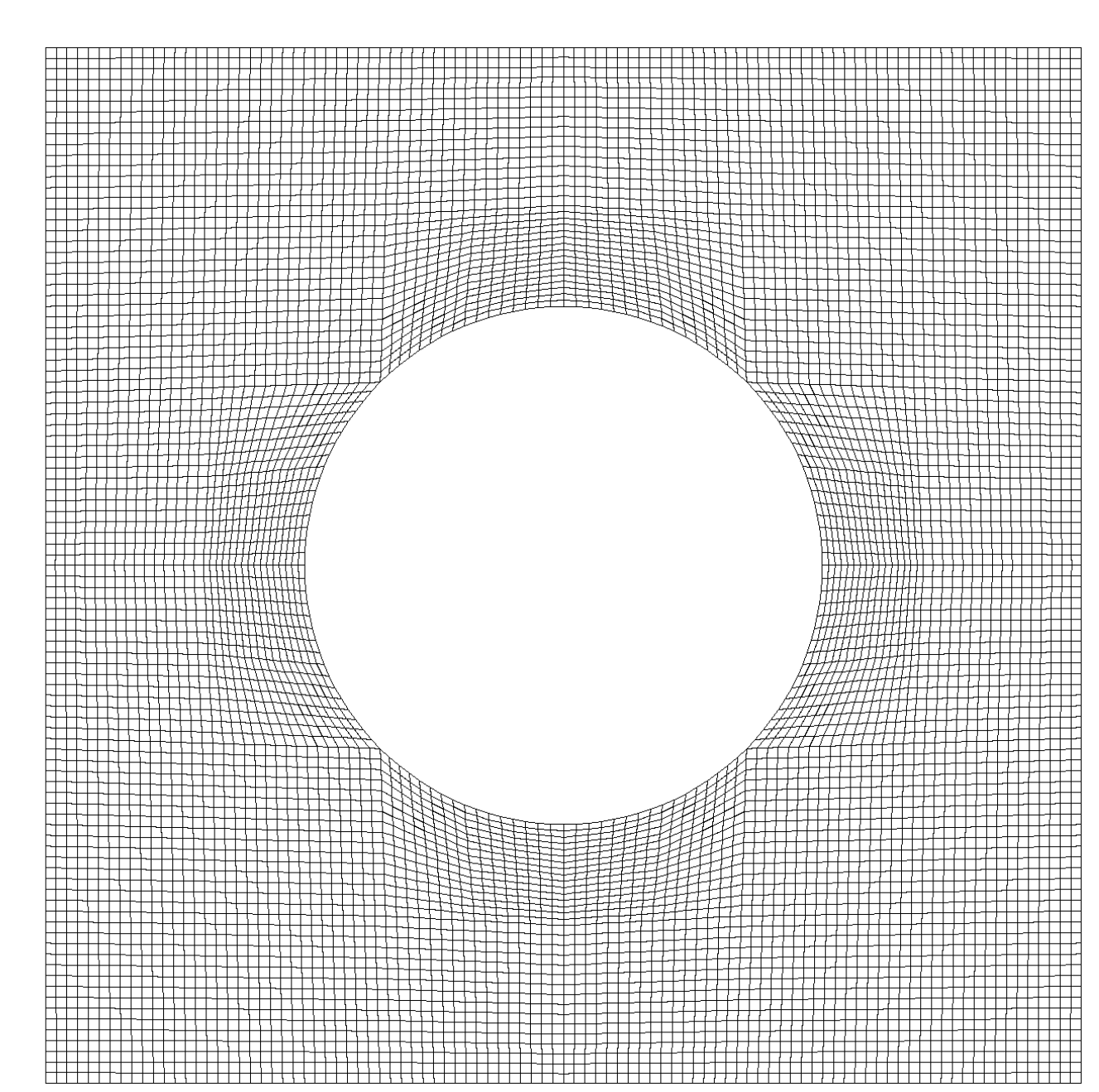}
\caption{Circle}\label{subfig.circle}
\end{subfigure}
\begin{subfigure}[htb]{0.3\textwidth}
\centering 
\includegraphics[trim=54mm 80mm 21mm 83mm, clip,width=0.8\textwidth]{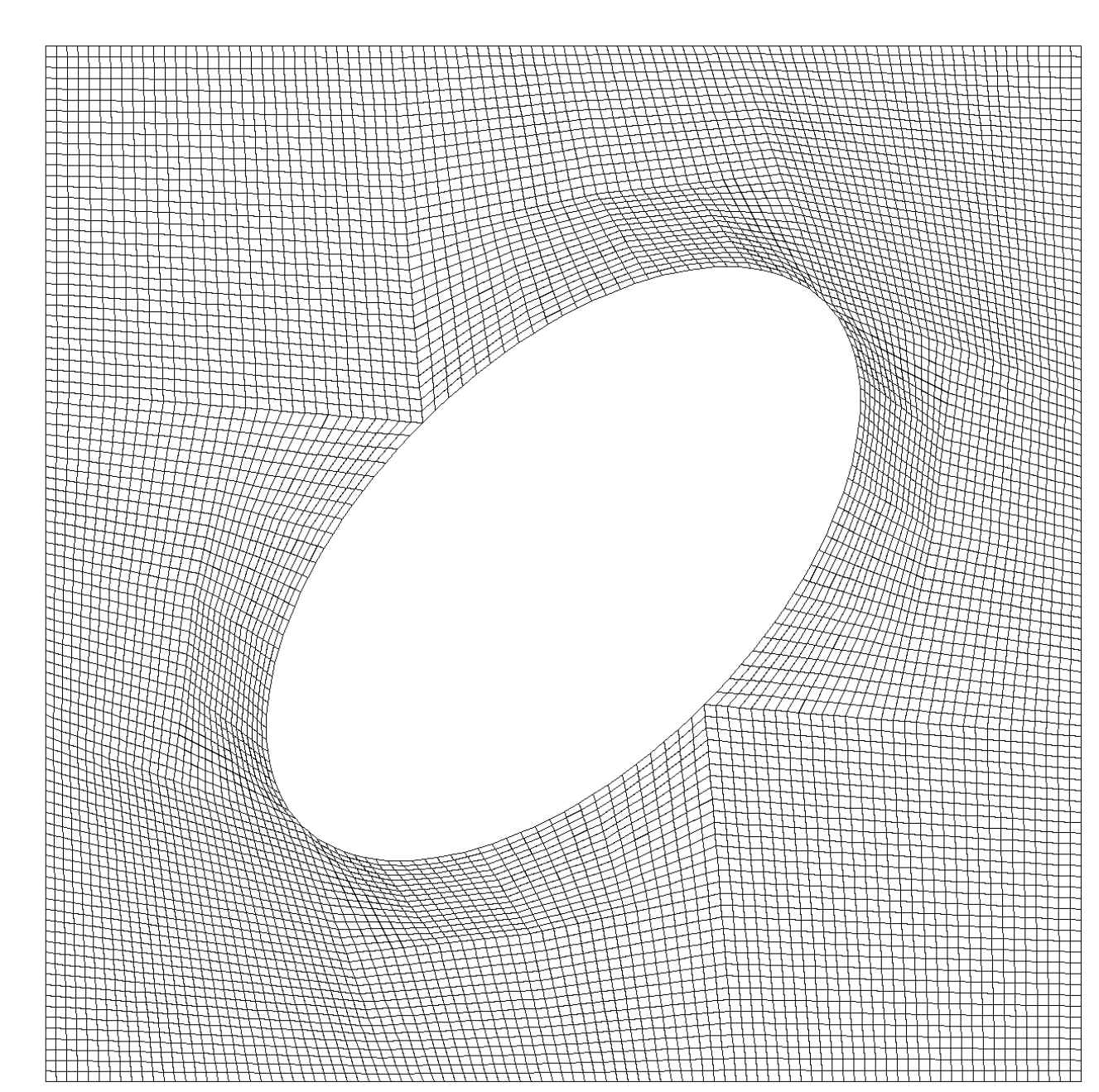}\caption{Ellipse}\label{subfig.ellipse}
\end{subfigure}
\caption{Mesh of the fluid part of the unit cell for the two types of inclusions: circles (\subref{subfig.circle}) and ellipses (\subref{subfig.ellipse}).}\label{fig.inclusions}
\end{figure}

As previously explained we use two different kinds of inclusion in the porous part, \textbf{circles} and \textbf{ellipses}. 
The geometries of the unit cells $Y = (0.1)^2$, see figure~\ref{fig.inclusions}, for these two cases are as follows:
\begin{enumerate}
\item the solid part of the cell $Y_s$ is formed by a circle with radius $0.25$ and center $(0.5, 0.5)$.
\item $Y_s$ consists of an ellipse with center $(0.5, 0.5)$ and semi-axes $a=0.357142857$ and $b=0.192307692$, which are rotated anti-clockwise by $45^\circ$. 
\end{enumerate}
The circle geometry is a case of axis symmetric geometry with respect to the axis $y$, perpendicular to the interface $\Sigma$, for which we expect from the theory that $\cblw=0$, see \cite*{JaegerMN:2001}.

All computations are done using the toolkit \texttt{DOpElib} (\cite{GollWW:2012}) based upon the C++-library \texttt{deal.II} (\cite{BangeHK:2007}).

\subsection{Finite element formulation of the Navier boundary layer problem}
\label{sec:Navier boundary problem}
The Navier boundary layer problem (\ref{BJ4.2})-(\ref{4.6}) is defined on $\Zbl:=Z^+ \cup \IF \cup Z^-$, where $\IF=(0,1)\times\{0\}$, $Z^+=(0,1)\times(0,+\infty)$ and $Z^-=\cup_{k=1}^\infty(Y_f - \{0,k\})$, with $Y_f$ the fluid part of the pore, see figure~\ref{subfig.unitcell}.

After \cite{JaegerM:1996} and \cite*{JaegerMN:2001}, it is known that $\boldsymbol\betabl$ converges exponentially towards $(\cbl,0)$ and $\wbl$ towards $\cblw$ in $Z^+$ for increasing $y_2$. On the porous side it has been also shown in the same references that the pressure $\wbl$ and the velocity $\boldsymbol\betabl$ converge exponentially towards zero. In addition it has been shown therein that
\begin{subequations}\label{eq.cbl_cblw}
\begin{align}
\cbl &= \int_0^1 \betabl_{1}(y_1,0) \ dy_1,\\
\cblw &= \int_0^1 \wbl(y_1,0) \ dy_1= \int_0^1 \wbl(y_1,a) \ dy_1,&&\forall a\geq 0,\label{equ.cbl_cblw:cblw}
\end{align}
\end{subequations}
where $(\boldsymbol\betabl, \wbl)$ is the solution of (\ref{BJ4.2})-(\ref{4.6}). Both this integrals are well defined since $\boldsymbol\betabl$ and $\wbl$ are smooth in $Z^+$ up to the interface $\Sigma$.

Since we can not deal with infinitely large domains, we consider a cut-off domain for numerical calculations defining the finite slab $\Zkl:=\Zbl\cap (0,1)\times(-l,k)$, $ k, l>0$. The distance of the cut-off from the interface, determined by $k$ and $l$, has to be taken large enough taking into account the exponential decay to reduce the approximation error introduced by cutting the domain. 

At the newly introduced parts of the boundary, namely $\Gamma_{ k} = (0,1)\times \{k\}$ and $\Gamma_{ l} = (0,1)\times \{-l\}$, we have to set some appropriate boundary conditions. We follow \cite*{JaegerMN:2001} and put zero Dirichlet condition for the two velocity components on $\Gamma_l$, while on $\Gamma_k$ a zero Dirichlet condition for the vertical component as well as zero normal flux of the first velocity component is imposed.


In the following we give the finite element approximation of the cut-off Navier boundary layer problem:
\begin{problem}[Cut-off Navier Boundary Layer]\label{prob.bl}
Find $\boldsymbol\beta^{bl}_h \in \tilde V_h(\Zkl)$ and $\omega^{bl}_h \in L_h(\Zkl)$, such  that 
 \begin{align}
 \int_{\Zbl} \left(\big(\nabla \boldsymbol\beta^{bl}_h + (\nabla \boldsymbol\beta^{bl}_h)^T\big) 
   \cdot \nabla \boldsymbol\varphi + \omega^{bl}_h \ \nabla \cdot \boldsymbol\varphi
   \right)\ dx&= -\int_{\IF}\varphi_1\ dx, && \forall \boldsymbol\varphi \in\tilde V_h(\Zkl),\\
 \int_{\Zbl} \nabla \cdot \boldsymbol\beta^{bl}_h \ \psi\ dx&= 0,  && \forall \psi \in L_h(\Zkl),
 \end{align}
\end{problem}
where the space for the velocity incorporates the Dirichlet boundary conditions on $\Gamma_k$ and $\Gamma_l$ and is thus defined as follows
\begin{align}
\tilde V_h(\Zbl) := \{\mathbf{v}_h \in C^0(\Zbl)~|~& {\mathbf v_h}_{|K} \in Q^2(K),K \in {\cal T}_h,\notag \\
&\mathbf v_h = (0,0) {\rm ~on} \cup_{n=1}^l (\partial Y_s - (0,n)),\notag \\
&\mathbf v_h = (0,0) \text{ on } \Gamma_l \text{ and } v_{h,1} =0 \text{ on } \Gamma_k,\notag\\
&\mathbf v_h {\rm ~is~} y_1-\mbox{periodic with period } 1 \}, \label{Space1}
\end{align}
with  $\partial Y_s$ the boundary of the inclusions in the pore domain, as shown in figure~\ref{subfig.unitcell}.

With the numerical solution of Problem~\ref{prob.bl} we approximate the
constants \cbl{} and \cblw{} for the considered inclusions. The
approximations $\cblh$ and $\cblwh$ are calculated by replacing in
\eqref{eq.cbl_cblw} the functions $(\boldsymbol\betabl, \wbl)$ with their discretized counterparts. As usual, the index $h$ indicates the approximation due to discretization. 
We  observe that to enhance the numerical approximation of \cblw{} in \eqref{equ.cbl_cblw:cblw} it is beneficial to calculate the integral of the pressure along a line far enough from the interface. We calculate the integral in \eqref{equ.cbl_cblw:cblw} for $a=1$, i.e. along the line $\Set{y\in \Zbl\ | \ y_2=1}$. 

The approximation of the cut-off problem by finite elements introduces two different sources of error: the \textbf{cut-off error} and the \textbf{discretization error}. In our computations, we set $k=l$ and compute the solution of Problem~\ref{prob.bl} for $1\leq k \leq 5$ on a family of hierarchic adaptively refined meshes.
To obtain the convergence results in Section~\ref{numerical_confirmation} with respect to $\epsilon$, it is important to control the cut-off and discretization errors and balance them to reduce the computational costs.
To balance the two errors we should cut the domain so that the order of the cut-off error equals that of the discretization error. To this aim we need to control the discretization error by a reliable estimation. Since we are interested on the calculation of \cbl{} and \cblw{}, we want to control directly the errors 
\begin{align*}
J_1=\cbl{} - \cblh, \quad J_\omega=\cblw{} - \cblwh.
\end{align*}
To this end, we employ the Dual Weighted Residual (DWR) method from \cite{BeckeR:2001} which gives an estimation of the discretization error with respect to a given functional (i.e. $J_1$ or $J_\omega$) exploiting the solution of a proper adjoint equation.
The DWR method in addition provides local error indicators to control the local mesh refinement.
The triangulation is adaptively refined until the estimated discretization error is smaller than a given tolerance. 

The reliability of this estimator has been shown in different applications in the context of flow problems and other problems, see e.g. \cite{Rannacher:99}, \cite{BeckeR:2001}, \cite{BraackR:2006}, \cite{Rannacher:2010}.
Nevertheless, we have performed an additional check to assure that the order of the error is indeed the one estimated.

To check the convergence we do not have the exact solution, but we can rely on the best approximation property of Galerkin approximations on quasi-uniform meshes to perform the following verification. 
On a series of uniformly refined grids we compute the approximations $C_{1,h}^{bl,unif}$ and $C_{\omega,h}^{bl,unif}$ and compare them with  reference values $C_{1,h}^{bl,ref}$ and $C_{\omega,h}^{bl,ref}$ computed on a (very fine) locally refined mesh. Additionally, we evaluate the error estimator $\eta$ on the uniformly refined grids and compare it with the following approximated errors: 
\[C_{1,h}^{bl,ref} - C_{1,h}^{bl,unif}, \quad  C_{\omega,h}^{bl,ref} -  C_{\omega,h}^{bl,unif}.\]
The results of this test show the expected reliability of the error estimator, since (cf. table~\ref{tab.results_global_mesh_refinement}) the solution on uniform meshes converges towards our reference solution and the error estimator is of the same order as the one given by the reference value. We have used grids with up to around 3.9 millions of degrees of freedom (DoF) for the verification with uniformly refined meshes.
Table~\ref{tab.results_global_mesh_refinement} shows the efficiency of the error estimator, i.e. 
\begin{align*}
I_{eff}(\cbl) = \frac{\eta(\cbl)}{C_{1,h}^{bl,ref} - C_{1,h}^{bl,unif}}, \quad I_{eff}(\cblw) = \frac{\eta(\cblw)}{C_{\omega,h}^{bl,ref} - C_{\omega,h}^{bl,unif}}.
\end{align*}
\begin{table}
\centering 
\begin{tabular}{r|ccc|ccc}
\toprule
\# DoF & $C_{1,h}^{bl,ref} - C_{1,h}^{bl,unif}$& $\eta_{\cbl}$ & $I_{eff}(\cbl)$& $C_{\omega,h}^{bl,ref} -  C_{\omega,h}^{bl,unif}$& $\eta_{\cblw}$ & $I_{eff}(\cblw)$\\
\cmidrule(lr){1-7}
	   1\,096 & -4.52E-04 & -2.54E-03 & 5.61 & 2.99E-02   &	-9.03E-03&	-0.30 \\
	   4\,142 & -7.49E-05 & -2.75E-04 & 3.67 & -2.07E-04  &  -1.02E-03 &	4.92  \\
	  16\,066 & -1.10E-05 & -1.87E-05 & 1.71 & 3.66E-06   &  -9.70E-05 &	-26.50\\
	  63\,242 & -9.60E-07 & -1.11E-06 & 1.15 & -9.72E-07  &  -4.78E-06 &	4.92  \\
	 250\,906 & -6.83E-08 & -7.60E-08 & 1.11 & -9.67E-08  &  -3.15E-07 &	3.26  \\
	 999\,482 & -4.54E-09 & -5.36E-09 & 1.18 & -9.03E-09  &  -2.77E-08 &	3.07  \\
	3\,989\,626 & -2.90E-10 & -3.82E-10 & 1.32 & -1.28E-09  &  -3.67E-09 &	2.85  \\
\bottomrule 
\end{tabular}
\caption{Results of the approximation of the constants \cbl{} and \cblw{} by uniform mesh refinement with $k=l=3$ and ellipses as inclusions. The first column gives the number of degrees of freedom (DoF).}\label{tab.results_global_mesh_refinement}
\end{table}
\begin{figure}
\centering
 \begin{subfigure}[b]{0.15\textwidth}
 \centering
\includegraphics[trim = 35 150 20 160,clip ,height=9cm]{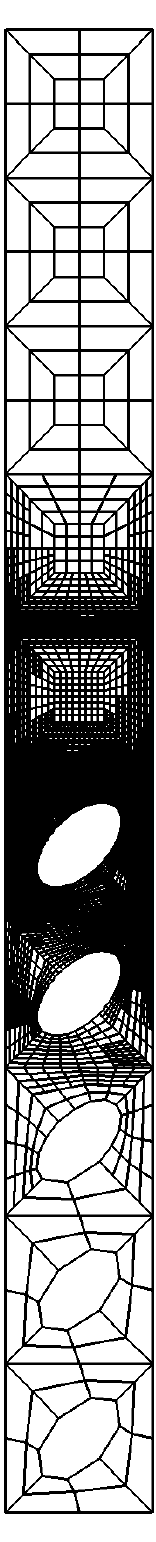}
\caption{Grid}
\label{fig.grid:overview}
\end{subfigure}
 \begin{subfigure}[b]{0.3\textwidth}
 \centering
\includegraphics[trim = 105 140 50 150,clip ,height=9cm]{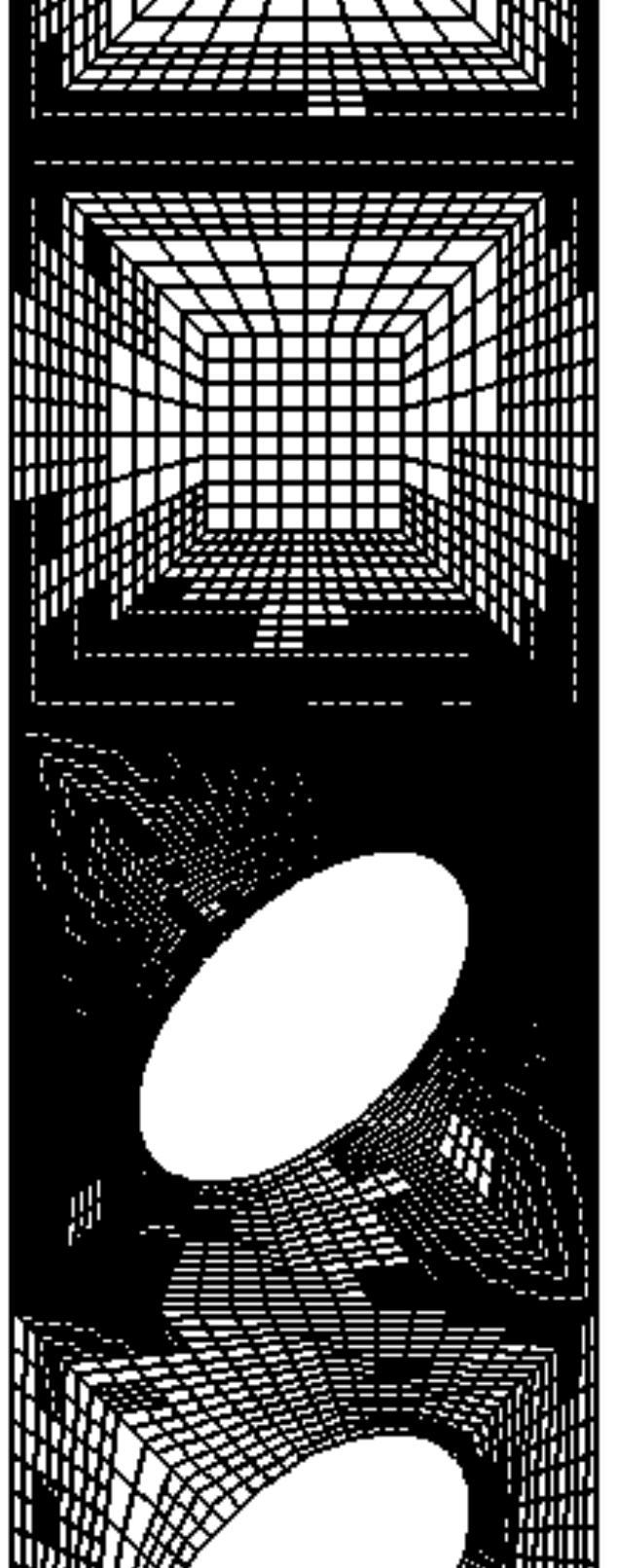}
\caption{Close-up}
\label{fig.grid:zoom}
\end{subfigure}
 \begin{subfigure}[b]{0.15\textwidth}
 \centering
\includegraphics[trim = 130 190 420 180,clip ,height=9cm]{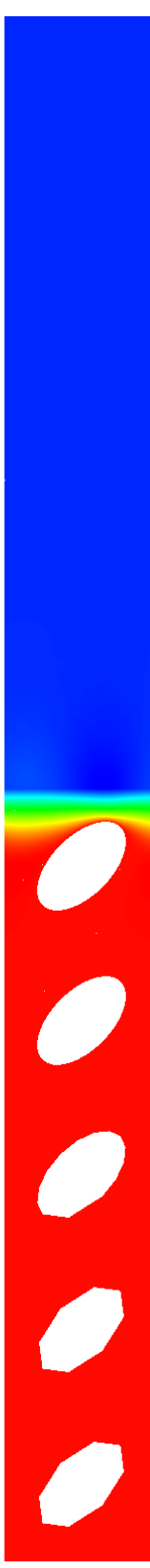}
\caption{$\beta_{h,1}$}
\label{fig.grid:v0}
\end{subfigure}
 \begin{subfigure}[b]{0.15\textwidth}
 \centering
\includegraphics[trim = 130 190 420 180,clip ,height=9cm]{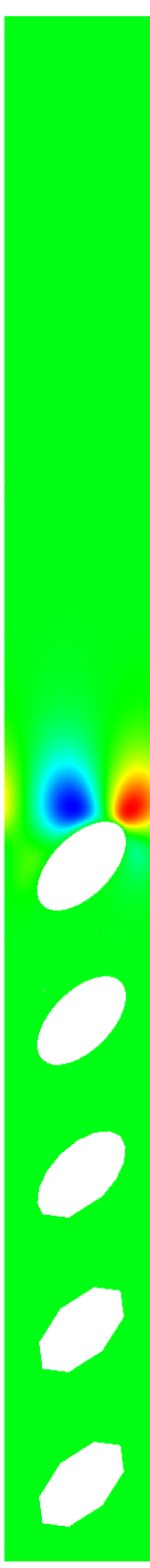}
\caption{$\beta_{h,2}$}
\label{fig.grid:v1}
\end{subfigure}
 \begin{subfigure}[b]{0.15\textwidth}
 \centering
\includegraphics[trim = 130 190 420 180,clip ,height=9cm]{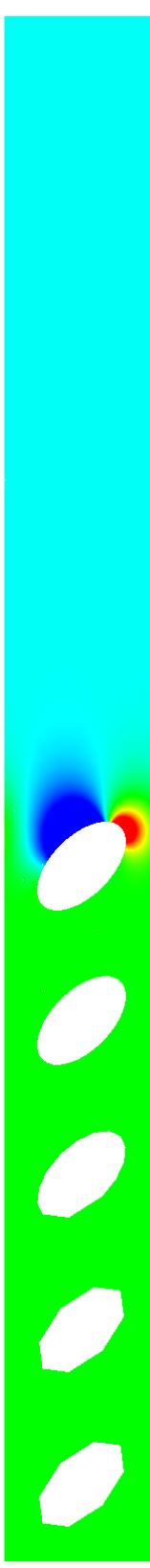}
\caption{$\omega_h$}
\label{fig.grid:pressure}
\end{subfigure}
\caption{Example of a locally refined grid for the adaptive computation of \cblw{} with $k=l=5$ in the Navier boundary layer problem. The whole mesh is shown in (\subref{fig.grid:overview}), whereas (\subref{fig.grid:zoom}) shows a close-up around the interface. In (\subref{fig.grid:v0}), (\subref{fig.grid:v1}) and (\subref{fig.grid:pressure}), the associated solution is shown.}\label{fig.grid}
\end{figure}
\begin{table}
\centering 
\begin{tabular}{r|cccc}
\toprule
k,l&\cblh& $\abs{\eta(\cbl{})}$ &\cblwh&$\abs{\eta(\cblw{})}$\\
\midrule
&\multicolumn{4}{c}{\textbf{circular inclusions}}\\
1  &-0.3038181652339	&	1.9E-12	&	-&-\\
2  &-0.3038219423526	&	2.0E-12	&	-&-\\
3  &-0.3038219423790	&	2.0E-12	&	-&-\\
4  &-0.3038219423789	&	2.0E-12	&	-&-\\
5  &-0.3038219423756    &	8.9E-13	&	-&-\\
\midrule
	&\multicolumn{4}{c}{oval inclusions}\\
1  & -0.2694539064491   & 4.3E-12	   &	-0.2413211012145	& 2.1E-11\\
2  & -0.2694545953967	& 1.4E-12      &	-0.2409146886571	& 7.3E-12\\
3  & -0.2694545953993	& 3.1E-12      &	-0.2409148310717	& 7.2E-12\\
4  & -0.2694545953993	& 2.1E-12      &	-0.2409148310975	& 8.5E-12\\
5  & -0.2694545953985	& 2.0E-12      &	-0.2409148310959	& 8.6E-12 \\ 		
\bottomrule 
\end{tabular}
\caption{Results of the approximation of the constants \cbl{} and \cblw{} as well as the estimated discretization error $\eta$ for different domain-lengths. }\label{tab.results_c1_cw}
\end{table}
Calculations of the two constants used in Section~\ref{numerical_confirmation} have been obtained setting the following tolerances $\eta(\cbl)$, $\eta(\cblw) < 10^{-11}$, where $\eta(\cbl)$ and $\eta(\cblw)$ are DWR error estimators respectively of $\cbl - C_{1,h}^{bl}$ and $\cblw - C_{\omega,h}^{bl}$. These tolerances are achieved by locally refined meshes with up to 7 millions of degrees of freedom. 

In figure~\ref{fig.grid:overview} an example of a mesh generated by the error estimator for the computation of \cblwh{} with $k=l=5$ is shown.  
A strong refinement can be observed in the neighborhood of the line $\Set{y\in\Zbl\ | \ y_2=1}$, where $\omega_h$ is evaluated to compute \cblwh{}, as well as in the vicinity of the first inclusion, see also figure~\ref{fig.grid:zoom} for a close-up of this region. 
In this part of the domain, the associated solution has large gradients, see figures~\ref{fig.grid:v0}, \ref{fig.grid:v1} and \ref{fig.grid:pressure} for an illustration  of the solution components. 
\begin{figure}
\centering
\resizebox{0.5\textwidth}{!}{\input{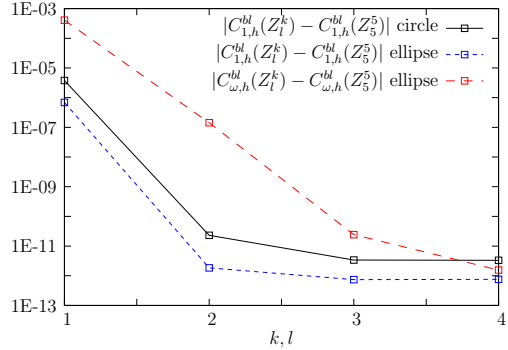}}
\caption{Difference between the computed constants on domains with increasing length and $\cblh(Z_5^5)$ resp. $\cblwh(Z_5^5)$.}\label{fig.exp_convergence}
\end{figure}
In \cite*{JaegerMN:2001} it is shown that the cut-off error decays exponentially with $k$ and $l$. 
To find the optimal cut-off level $l, k$ we perform a convergence check taking as reference value the constants computed on $Z_5^5$, i.e. $\cblh(Z_5^5)$ resp. $\cblwh(Z_5^5)$.
Figure~\ref{fig.exp_convergence} shows the error between the constants computed on $\Zkl$ with $k=l =1,\dots, 4$ and the reference values computed on $Z_5^5$. 
The exponential decay of the cut-off error with the distance from interface can be observed for both approximations. Furthermore, it can be observed that the error $|\cblh(Z_k^l) - \cblh(Z_5^5)|$ is of the order of the discretization error, i.e. $10^{-12}$, for $k,l\geq3$ for \cbl{} and $k,l\geq 4$ for \cblw{}.

In the following, we use approximations computed with local mesh refinement and $k,l=5$ maintaining the simplified notation \cblh{} and \cblwh{}.
For details on the used local refinement strategy see \cite{Richter:Diss}.
The calculated constants and the respective error estimation are listed in table~\ref{tab.results_c1_cw}.
\subsection{Cell problem and determination of the permeability}
\label{subsec.cell problem}
For a numerical confirmation of 
\eqref{ConcPression} we need the solution of appropriate cell problems, depending on the shape of the inclusions, to calculate the rescaled permeability $K$.

To introduce the weak form of the cell problems, we define the following function space
\begin{align}
\hat V(Y_f) &:= \Set{\mathbf{v} \in H^1(Y_f)^2\ | \ \mathbf{v}  =0 \; \hbox{on } \ \p Y_s, \mathbf{v} \hbox{ is } 1\text{-periodic}}.
\end{align}
Following the derivation of Darcy's law by homogenization, the matrix $K$ is defined as
\begin{align}
K_{ij} = \int_{Y_f} w_j^i\ dy,\quad i,j=1,2,
\end{align}
where $\mathbf w$ is the velocity of the following \textbf{cell problem}.
\begin{problem}[Cell Problem]\label{prob.cell_problem}
Let $i,j=1,2$. Find a velocity field $\mathbf w^i \in \hat V(Y_f)^2$ and a pressure $\pi^i \in L_0(Y_f)$, such that,
\begin{align}
\int_{Y_f} \left(\big(\nabla \mathbf w^{i}+(\nabla \mathbf w^{i})^T \big)\cdot \nabla \boldsymbol\varphi + \pi^{i} \ \nabla \cdot \boldsymbol\varphi\right)\  dx&= \int_{Y_f}\varphi_i\ dx, && \forall \boldsymbol\varphi \in \hat V(Y_f),\\
 \int_{Y_f} \nabla \cdot \mathbf w^i \ \psi\ dx&= 0,  && \forall \psi \in L_0(\Zbl).
\end{align}
\end{problem}

The cell problem is solved with Taylor-Hood elements and an adaptive algorithm based on the DWR method to compute precisely the reference values for the permeability matrix $K$, see also Subsection~\ref{sec:Navier boundary problem}. Each component $w^i$ is solved by a tailored grid refinement considering as goal functional for the a posteriori error estimation the components of $K$.

The computed reference values for the circles are
\begin{align}  K^{circ}_h &= k^{circ}_h Id \approx {0.01990143534975}Id
\end{align}
with an estimated discretization error of $1.38 \ 10^{-11}$. For the case with ellipses as inclusion the following values have been calculated
\begin{align}
K^{oval}_h &= \left(
\begin{array}{cc}
K_{h,11} & K_{h,12} \\
K_{h,12} & K_{h,22} \\
\end{array}
\right) \approx \left(
\begin{array}{cc}
0.0159787174788 &  0.00303449804138  \\
0.00303449804138 & 0.0159787174788
\end{array}\right).
\end{align}
The estimated discretization errors are $2.76 \ 10^{-12}$ for $K_{h,11}$ and $1.10 \ 10^{-13}$ for $K_{h,12}$.

In the next session we use the reference values of \cblh{}, \cblwh{} and $K_h$ to present a numerical confirmation of the interface law.

\section{Numerical confirmation of the interface law}
\label{numerical_confirmation}
In the context of the coupling between Stokes and Darcy, the slip condition for the velocity of the free flow \eqref{4.95} has been established numerically for example in \cite{Kaviany:1995, SahraK:1992} and \cite{LarsoH:1986, LarsoH:1987}. In addition, the numerical calculation of the constants $C^{bl}_1$ and $C^{bl}_\omega$ has been performed by finite element method in \cite{JaegerMN:2001}. Nevertheless, numerical results on the evidence of the pressure relation \eqref{Presspm2A} based on a comparison between the microscopic and the homogenized flow has not yet been shown. This is the goal of this section.
Note that \cite{Kaviany:1995, SahraK:1992} always deal with isotropic geometries and, consequently, for small Reynolds numbers they do not observe pressure jump, see also Section~\ref{intro}.

Specifically we show numerical evidence based on the following consideration.
As will be clear from the results of this section, in the microscopic model the pressure values,  which converge to $p_D$ on $\O_p$ and to $p^\eff$ on $\O_f$, oscillate due to the inclusions. 
Approaching the interface, the microscopic pressure oscillates with amplitude which does not vanish with $\epsilon$.
Indeed, it has been shown in \cite{MarciM:2012} that the pressure of the microscopic model on the interface converges to $p_\eff$ in the sense of bounded measures that allows such oscillations. We give a numerical justification of the interface laws for the pressure and shear stress averages over the pore faces at the interface as explained in more detail in Subsection~\ref{sec:case_I}.

We first focus on the theoretical results (\ref{4.100A}--\ref{4.100}) and
(\ref{ConcPression}), which take in consideration quantities defined over
the domains $\O_p$ or $\O_f$ and can be used to show the estimates  \eqref{equ.bjs epsilon32} and \eqref{EstPressBdry}.

At the interface the estimates \eqref{equ.bjs epsilon32} and \eqref{EstPressBdry} describe precision of the Beavers and Joseph slip condition and order of approximation for the pressure field. We refer to
\cite{MarciM:2012} for details on these estimates.

In next subsections, we consider two different flow conditions: a periodic
flow and the flow with an injection condition, called here Beavers-Joseph
case.

\subsection{Case I: periodic case}
\label{sec:case_I}
In this subsection we present the \textit{periodic case}, i.e. boundary conditions and microstructure of the porous domain are periodic.

To show convergence with epsilon without the predominance of the discretization errors, we compute the solutions of the microscopic Problem~\ref{prob.micro_weak} for  
\begin{align*}
\epsilon\in \Set{1,\frac 1 3, 0.1, \frac{1}{31}, 0.01,\frac 1{316}, 0.001, \frac1 {3162}}.
\end{align*} 
The data of the simulations are $\mathbf f=(1,0)$, $L = 1$, and the rigid inclusions of the porous part are either ellipses or circles with geometry described in Section~\ref{numerical_setting}, see also figure~\ref{fig.inclusions}.

In the periodic case we can exploit the fact that the solutions $u^\ep$ are not only $L$-periodic but, due to the constant right hand side, also $\epsilon$-periodic in $x_1$-direction. To strongly reduce computational costs we compute the approximations $\mathbf{\tilde u}_h^\ep$ on a domain with length $\epsilon$ instead of $L$. We employ then the $\epsilon$-periodicity to reconstruct the solution $\mathbf u_h^\ep$ on the whole domain with length $L$.
The computational grids are obtained by global refinement.

To show the convergence with $\epsilon$ in (\ref{4.100A}--\ref{4.100}) and (\ref{ConcPression}), we need the solution of the microscopic problem, the constants \cbl{} and \cblw{}, and the permeability $K$. In addition, we can use the expression of the exact solution of the effective/macroscopic flow: $u^{\eff}$, $p^{\eff}$ and $p_D$.
The analytical solution for the effective problem (\ref{4.91})-(\ref{4.95}) is
\begin{gather}
u_1^{\eff} (x) = \frac{1-x_2}{2} \frac{x_2 (1-\ep C^{bl}_1 )- \ep C^{bl}_1 }{1-\ep C^{bl}_1} ,\label{horvel} \\
  u_2^{\eff} =0 \quad \mbox{and} \quad p^{\eff} =0. \label{varothers}
\end{gather}
It follows
\begin{align}
\sigma^{\eff}_{12} (0) = \frac{\partial u_1^{\eff}}{\partial x_2} (0) = \frac{1}{2  (1-\ep C^{bl}_1 )}\quad \text{ and } \quad M^{\eff} = \frac{1}{12} \frac{1-4\ep C^{bl}_1 }{1-\ep C^{bl}_1}.
\end{align}
We have then in the porous medium $\O_p$
\begin{equation}\label{porpressure}
    p_D (x) = \frac{C^{bl}_\o}{2  (1-\ep C^{bl}_1 )} + \frac{K_{12}}{K_{22}} x_2 = \frac{C^{bl}_\o}{2}+ \frac{K_{12}}{K_{22}} x_2 + O(\ep).
\end{equation}
The analytical solutions are evaluated using the approximated values \cblh{}, \cblwh{} as well as the permeability $K_h$ computed in the previous Section~\ref{numerical_setting} with a discretization error of the order at least $10^{-12}$. 
A discretization error is thus included in the calculation of $u_1^{\eff}$ and $p_D$, but it is negligible in comparison with the error in $\epsilon$ for the values considered in our convergence tests.
To ease notation, we do not distinguish between the analytical solutions and their approximations due to discretization errors of the constants used in the expressions.

We have now everything at hand to compute the error estimates between the solutions of the microscopic problems and the solutions of the effective/macroscopic  problems.  Direct simulations confirm (\ref{4.100A})-(\ref{4.100}), i.e.
\begin{gather}\label{eq.v_and_m_con_fluid}
\int_{\O_f} | \ve - u_1^{\eff} (x_2 ) \mathbf{e}^1 |^2 \ dx   + \vert  M^\ep - \frac{1}{12} \frac{1-4\ep C^{bl}_1 }{1-\ep C^{bl}_1}  \vert^2 = O (
\ep^{3}),  \\
\int_{\O_f} \{  |  p^\ep  |
+ |  \nabla \ve -  \bigg( \frac{1}{2  (1-\ep C^{bl}_1 )} - x_2 \bigg)  \left(
                                                    \begin{array}{cc}
                                                      0 & 1 \\
                                                      0 & 0 \\
                                                    \end{array}
                                                  \right)
 |\} \ dx = O (\ep).\label{eq.p_and_nabla_v_con_fluid}
\end{gather} 
as can be seen in {table~\ref{tab.eps_con_domains} or {figure~\ref{fig.eps_con_domains}}. 

In addition, in the porous medium direct simulations confirm (\ref{ConcPression}) 
 \begin{align} \label{equ.porous_pressure_estimate}
 \int_{\O_p} | p^\ep - \frac{C^{bl}_\o}{2} - \frac{K_{12}}{K_{22}} x_2 |^2 \ dx = o(1),
 \end{align}
see the last column in table~\ref{tab.eps_con_domains} as well as figure~\ref{fig.eps_con_domains.1storder_circle} and \ref{fig.eps_con_domains.1storder_ellipse}. In the periodic case convergence of order $O(\epsilon)$ can be observed.

\begin{figure}
\centering
 \begin{subfigure}[htb]{0.48\textwidth}
 \centering
 \resizebox{0.9\textwidth}{!}{\input{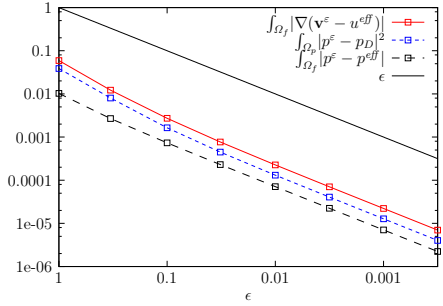}}
 \caption{Confirmation of \eqref{eq.p_and_nabla_v_con_fluid} and \eqref{equ.porous_pressure_estimate}, circles as inclusions.}
   \label{fig.eps_con_domains.1storder_circle}
 \end{subfigure}%
 \quad
 \begin{subfigure}[htb]{0.48\textwidth}
 \centering
 \resizebox{0.9\textwidth}{!}{\input{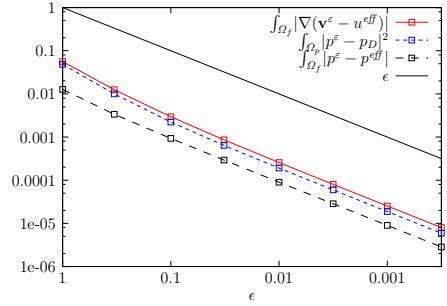}}
 \caption{Confirmation of \eqref{eq.p_and_nabla_v_con_fluid} and \eqref{equ.porous_pressure_estimate}, ovals as inclusions.}
   \label{fig.eps_con_domains.1storder_ellipse}
 \end{subfigure}
 
  \begin{subfigure}[htb]{0.48\textwidth}
 \centering
 \resizebox{0.9\textwidth}{!}{\input{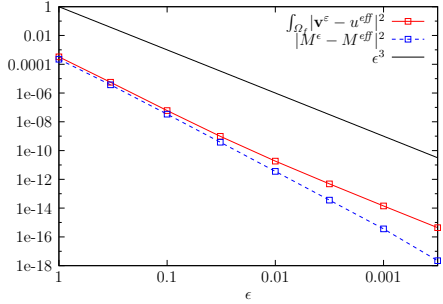}}
 \caption{Confirmation of \eqref{eq.v_and_m_con_fluid}, circles as inclusions.}
 \end{subfigure}%
 \quad
 \begin{subfigure}[htb]{0.48\textwidth}
 \centering
 \resizebox{0.9\textwidth}{!}{\input{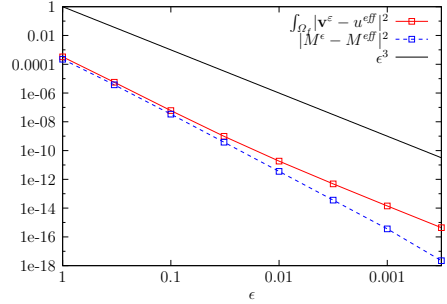}}
 \caption{Confirmation of \eqref{eq.v_and_m_con_fluid}, ovals as inclusions.}
 \end{subfigure}
 \caption{{Confirmation of the estimates \eqref{eq.v_and_m_con_fluid}, \eqref{eq.p_and_nabla_v_con_fluid} and \eqref{ConcPression} for oval (right column) and circular (left column) inclusions. Notice the different logarithmic scaling in the two rows.}}
 \label{fig.eps_con_domains}
 \end{figure}
{
 \begin{table}
\centering 
\begin{tabular}{r|ccccc}
\toprule
{$\epsilon$}	&\multicolumn{1}{c}{$\int_{\O_f}\abs{\ve - \mathbf u^{\eff}}^2$} &\multicolumn{1}{c}{$\abs{M^\epsilon - M^\eff}$}&\multicolumn{1}{c}{$\int_{\O_f}\abs{p^\ep-p^\eff}$}	& \multicolumn{1}{c}{$\int_{\O_f}\abs{\nabla (\ve-\mathbf u^\eff)} $} & \multicolumn{1}{c}{$\int_{\O_p}\abs{p^\ep-p_D}^2$}\\
\midrule
&\multicolumn{5}{c}{\textbf{circular inclusions}}\\
1      		&	4.69E-04 	&	1.83E-02 	&	1.03E-02 	&	5.92E-02	&	3.87E-02 \\
$\frac 1 3$	& 	8.22E-06	& 2.40E-03 		&2.69E-03		&1.22E-02	    & 8.05E-03\\
0.1    		&	8.28E-08 	&	2.31E-04 	&	7.36E-04 	&	2.72E-03	&	1.64E-03 \\
$\frac 1{31}$ &1.17E-09		& 2.46E-05		&2.31E-04		&7.69E-04		&4.50E-04\\
0.01   		&	1.81E-11 	&	2.38E-06 	&	7.09E-05 	&	2.26E-04	&	1.31E-04 \\
$\frac 1 {316}$&4.08E-13	&2.38E-07		&2.24E-05		&7.04E-05		&4.06E-05\\
0.001  		&	1.12E-14 	&	2.38E-08 	&	7.06E-06 	&	2.21E-05	&	1.28E-05 \\
$\frac 1 {3162}$&3.36E-16	&1.79E-09		&2.26E-06		&7.00E-06		&4.03E-06\\
\midrule
	&\multicolumn{5}{c}{\textbf{oval inclusions}}\\
1 				&	3.25E-04	&	1.49E-02 	&	1.28E-02  	&  	5.61E-02 & 4.83E-02 \\
$\frac 1 3$		&5.64E-06		&1.93E-03 		&3.38E-03		&1.25E-02	 &1.00E-02\\
0.1				&	6.03E-08	&	1.85E-04 	&	9.33E-04  	&  	2.96E-03 & 2.23E-03 \\
$\frac 1{31}$ 	&9.83E-10		&1.95E-05		&2.93E-04		&8.64E-04	 &6.42E-04\\
0.01 			&	1.86E-11 	&	1.89E-06 	&	9.02E-05  	&  	2.58E-04 & 1.91E-04 \\
$\frac 1 {316}$	&4.84E-13		&1.90E-07 		&2.85E-05		&8.05E-05	 &5.98E-05\\
0.001			&	1.42E-14 	&	1.89E-08 	&	8.99E-06  	&  	2.53E-05 & 1.89E-05 \\
$\frac 1 {3162}$& 4.37E-16		&1.07E-09		&2.85E-06		&8.00E-06	 &6.03E-06\\
\bottomrule 
\end{tabular}
\caption{Confirmation of the estimates \eqref{eq.v_and_m_con_fluid},
\eqref{eq.p_and_nabla_v_con_fluid} and \eqref{equ.porous_pressure_estimate} for oval and circular inclusions.}\label{tab.eps_con_domains}
\end{table}

\begin{figure}
\centering
\begin{subfigure}[htb]{0.45\textwidth}   
\includegraphics[trim = 80mm 80mm 80mm 80mm,clip, width=\textwidth]{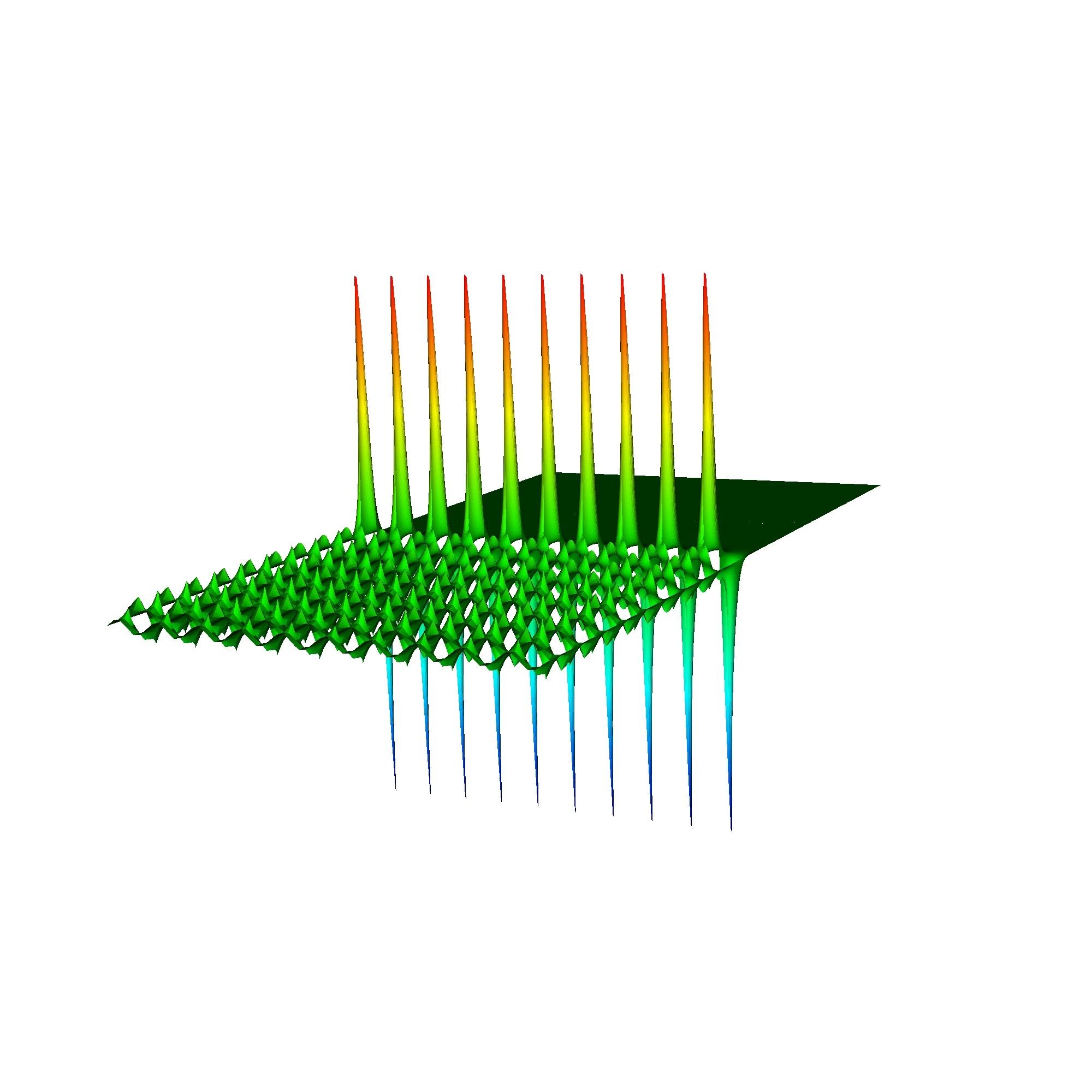}
\caption{Circles as inclusions.}\label{fig.pressure_per_em1.circle}
\end{subfigure}
\begin{subfigure}[htb]{0.45\textwidth}
\centering
 \includegraphics[trim = 80mm 80mm 80mm 80mm, clip, width=\textwidth]{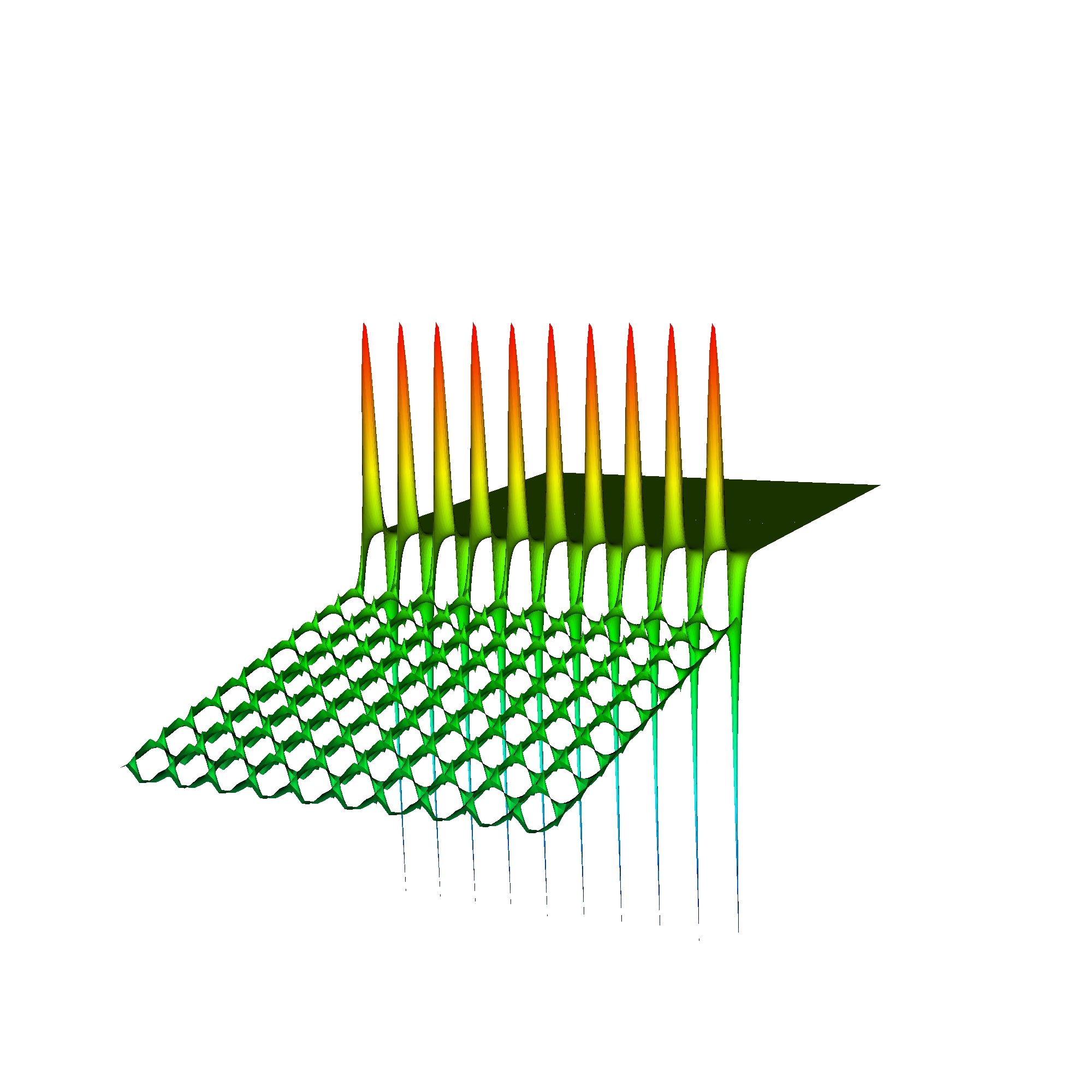}
 
 \caption{Ellipses as inclusions.}\label{fig.pressure_per_em1.ellipse}
 \end{subfigure}
\caption{Visualization of the pressure $p^\ep$ with $\epsilon = 10^{-1}$ and periodic boundary conditions.}
\label{fig.pressure_per_em1}
\end{figure}

\begin{remark}[Pressure peaks]
\label{remark:pressure_peaks}
The pressure in the porous domain oscillates due to the inclusions. 
In particular, we see in figure~\ref{fig.pressure_per_em1} that the pressure at the boundary of each inclusion adjacent to the interface has two prominent peaks, one positive and one negative.
To visualize it better, we refer to figure~\ref{fig.pressure_ellipse_per_em0_zoom}, where we show the pressure for $\epsilon = 1$. It can be clearly seen that the solution is smooth and bounded and the maximum and minimum values are on the boundary, as expected by the maximum principle.
\end{remark}

\begin{figure}
\centering
 \begin{subfigure}[htb]{0.48\textwidth}
 \centering
 \includegraphics[trim = 80mm 80mm 80mm 80mm,clip, width=0.8\textwidth]{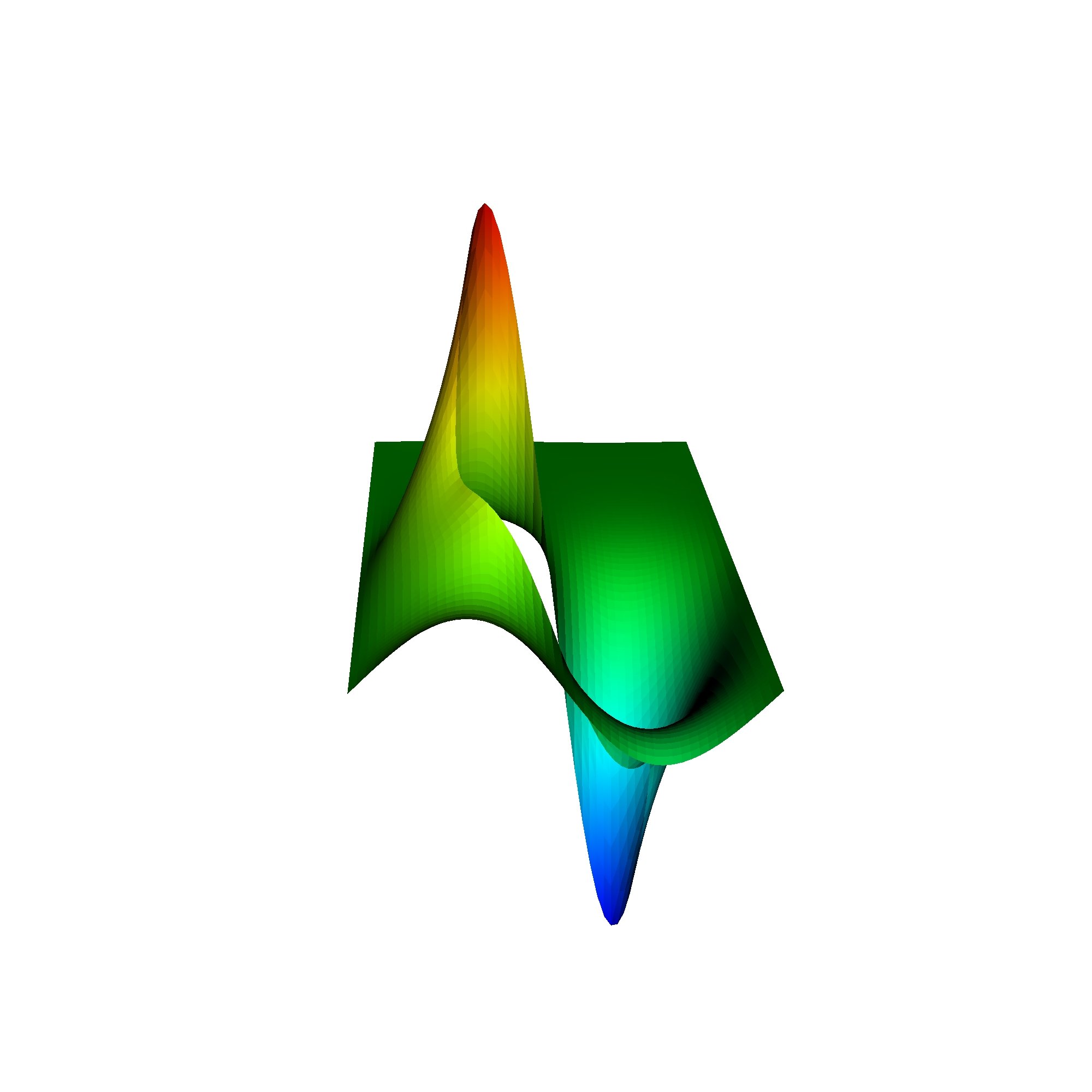}
 \caption{Circles as inclusions.}
 \end{subfigure}%
  \begin{subfigure}[htb]{0.48\textwidth}
 \centering
 \includegraphics[trim = 80mm 80mm 80mm 80mm, clip, width=0.8\textwidth]{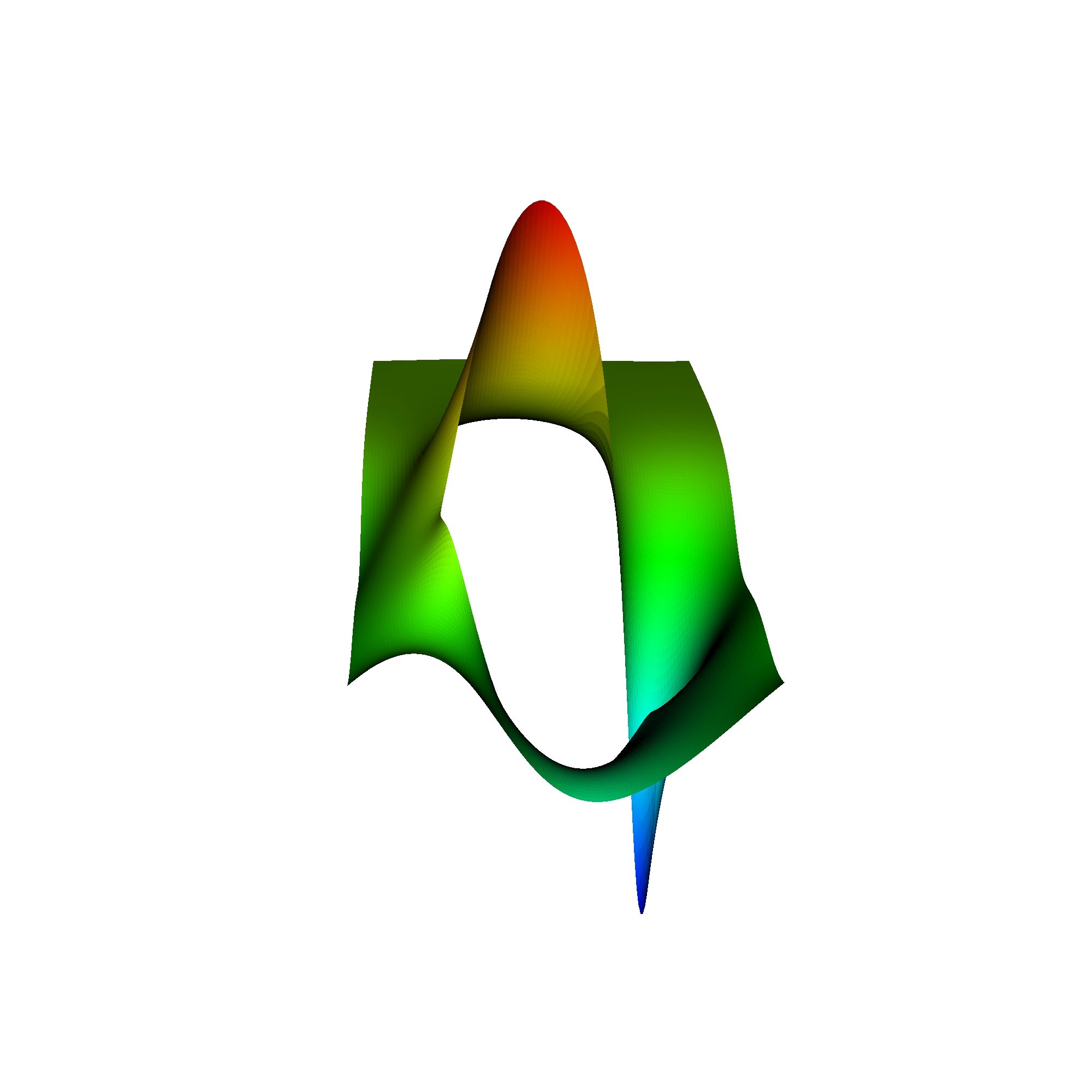}
 \caption{Ellipses as inclusions.}
 \end{subfigure}%

\caption{Plot of the pressure values with $\epsilon = 1$ and periodic boundary conditions.}
\label{fig.pressure_ellipse_per_em0_zoom}
\end{figure}

\begin{figure}
\centering
\resizebox{0.6\textwidth}{!}{\input{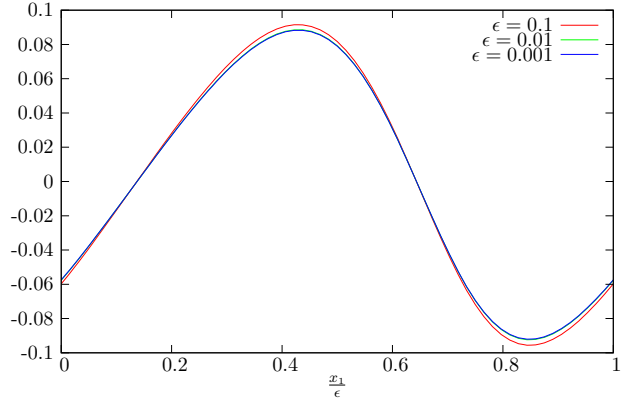}}
\caption{Values of the pressure $p^\ep$ on the interface $\Gamma$ for different $\epsilon$. 
The horizontal axis is scaled by $\epsilon$, i.e. the plot shows $p^\ep(x_1/\epsilon,0)$ for $x_1\in[0,\epsilon]$.}
\label{fig.p_on_gamma}
\end{figure}


As can be observed in figure~\ref{fig.pressure_per_em1}, in case of oval inclusions, in the vicinity of the interface appears a pressure jump. 
In the periodic case the interface law is rigorously confirmed by the convergence
rates of table~\ref{tab.eps_con_domains}. Nevertheless, in the following we
give an insightful illustration of the jump behavior on the interface. This
procedure is tested in the periodic case, which is supported by theoretical
results, and in the next subsection is applied to a more general flow
condition.

The amplitude of the pressure oscillations on the interface is of order $O(1)$ with respect to $\epsilon$. Figure~\ref{fig.p_on_gamma} clearly shows this behavior. Note that $p^\epsilon$ is depicted for three different $\epsilon$-domains scaled with $\epsilon$ for comparison purposes.
Furthermore, we have observed that the average value of $p^\epsilon$ over one period converges towards zero, which is the value of $p^\eff$.
As a consequence of this observation, to define the jump we introduce cell-wise averaged quantities on $\Gamma$ denoted with a bar over it $\clos{\cdot}$. Let $x \in \Gamma$, $m\in \mathbb N$ with $m\epsilon<x<(m+1)\epsilon$, the cell average of $\partial_2u^\ep_1$ is then defined as
\begin{align}
\clos{\partial_2u^\ep_1}(x) = \frac 1 \epsilon \int_{m\epsilon}^{(m+1)\epsilon}\partial_2 u^\ep_1(s, 0)\ ds.
\end{align}
The values $\clos u^\ep_1$ and $\clos p^\ep$ are defined analogously.
The pressure jump for the continuous microscopic pressure is defined by values taken on $\Gamma$ and on a line below $\Gamma$. The distance of this line from the interface is heuristically motivated by the following consideration. 
In view of the continuity of $p^\ep$, the ``pressure from below'' has to be taken away from the interface. In addition, since the pressure in $\Omegap$ converges to an affine function (cf. \ref{porpressure}), the line below must be not too far from the interface to define the jump.
We define hence ``the pressure from below'' $\clos p^\ep_{\operatorname{d}}$ by 
\begin{align}
\clos p^\ep_{\operatorname{d}}(x)= \frac 1 \epsilon \int_{m\epsilon}^{(m+1)\epsilon}p^\ep(s, -2\epsilon)\ ds.
\end{align} 
As a remark, calculations have been done for the pressure taken at a
distance $\epsilon$, $2 \epsilon$ and $3 \epsilon$ from the interface in
the porous part. In these three tests we observed convergence with epsilon, but only the distance $2 \epsilon$, in this specific case, gives the perfect convergence rate as can be observed in table~\ref{tab.eps_con_mean}, where
\begin{align}
\norm{g}_{2,\Gamma}=\sqrt{\int_\Gamma\abs{g(x)}^2\ dx}.
\end{align}
In table~\ref{tab.eps_con_mean} the third and fourth columns show a
convergence of order $O(\epsilon)$ for the Beavers-Joseph and jump
interface condition.

We use this heuristic definition of the jump also in the next subsection for a more general case.

\begin{table}
\centering
\begin{tabular}{r|cccc}
\toprule
$\epsilon$	&$\norm{\clos{\partial_2u^\ep_1}-\frac{1}{2(1-\ep \cbl)}}_{2,\Gamma}$&$ \norm{\clos{p^\ep_{\operatorname{d}}}-\frac \cblw 2}_{2,\Gamma}$&$\norm{\clos{p^\ep_{\operatorname{d}}}-\clos{p^\ep}-\cblw\clos{\partial_2u^\ep_1}}_{2,\Gamma}$&$ \norm{\frac{\clos{u_1^\ep}}{\ep}+\cbl\clos{\partial_2u^\ep_1}}_{2,\Gamma}$\\
\midrule
&\multicolumn{4}{c}{\textbf{circular inclusions}}\\
   1             	& 3.66E-02 &          -  &  -  		  &  4.77E-02\\
  $\frac 13 $     	& 4.81E-03 &   7.95E-14  &  7.20E-14  &  1.59E-02\\
   0.1           	& 4.63E-04 &   5.70E-14  &  5.68E-14  &  4.77E-03\\
 $ \frac 1 {31}$  	& 4.91E-05 &   9.88E-14  &  8.66E-14  &  1.54E-03\\
   0.01          	& 4.75E-06 &   4.96E-14  &  4.06E-14  &  4.77E-04\\
 $ \frac 1 {316}$ 	& 4.77E-07 &   8.66E-14  &  7.95E-14  &  1.51E-04\\
   0.001         	& 4.76E-08 &   7.72E-14  &  6.94E-14  &  4.77E-05\\
 $ \frac1 {3162} $	& 6.68E-09 &   3.02E-14  &  2.77E-14  &  1.51E-05\\
\midrule
&\multicolumn{4}{c}{\textbf{oval inclusions}}\\
   1               	&  2.99E-02  &    -         &   -         &  3.79E-02\\
  $ \frac 13  $     &  3.86E-03  &    1.29E-01  &   1.39E-01  &  1.26E-02\\
   0.1             	&  3.69E-04  &    3.86E-02  &   4.18E-02  &  3.79E-03\\
 $  \frac 1 {31} $  &  3.91E-05  &    1.24E-02  &   1.35E-02  &  1.22E-03\\
   0.01            	&  3.78E-06  &    3.86E-03  &   4.18E-03  &  3.79E-04\\
 $  \frac 1 {316} $	&  3.79E-07  &    1.22E-03  &   1.32E-03  &  1.20E-04\\
   0.001           	&  3.78E-08  &    3.86E-04  &   4.19E-04  &  3.79E-05\\
 $ \frac 1 {3162} $	&  5.68E-09  &    1.22E-04  &   1.32E-04  &  1.20E-05\\
\bottomrule 
\end{tabular}
 \caption{Numerical evidence of the convergence results on the interface $\Gamma$ in the case of periodic boundary conditions. }\label{tab.eps_con_mean}
\end{table}

\subsection{Case II: Beavers-Joseph case}
\label{sec:case_II}
In this case we investigate the behavior of the microscopic solutions for a set of boundary conditions corresponding to the original experiment of \cite{BeaveJ:1967}. 
As a difference from the periodic case, a pressure drop in $x_1$-direction is prescribed and the vertical velocity component is set to zero on $\Gamma_\per$. On the same boundaries, due to the divergence free condition, it holds $\partial_1 u_1 = 0$.
\begin{figure}
\centering
\begin{subfigure}[b]{0.48\textwidth}
\includegraphics[width=0.98\textwidth]{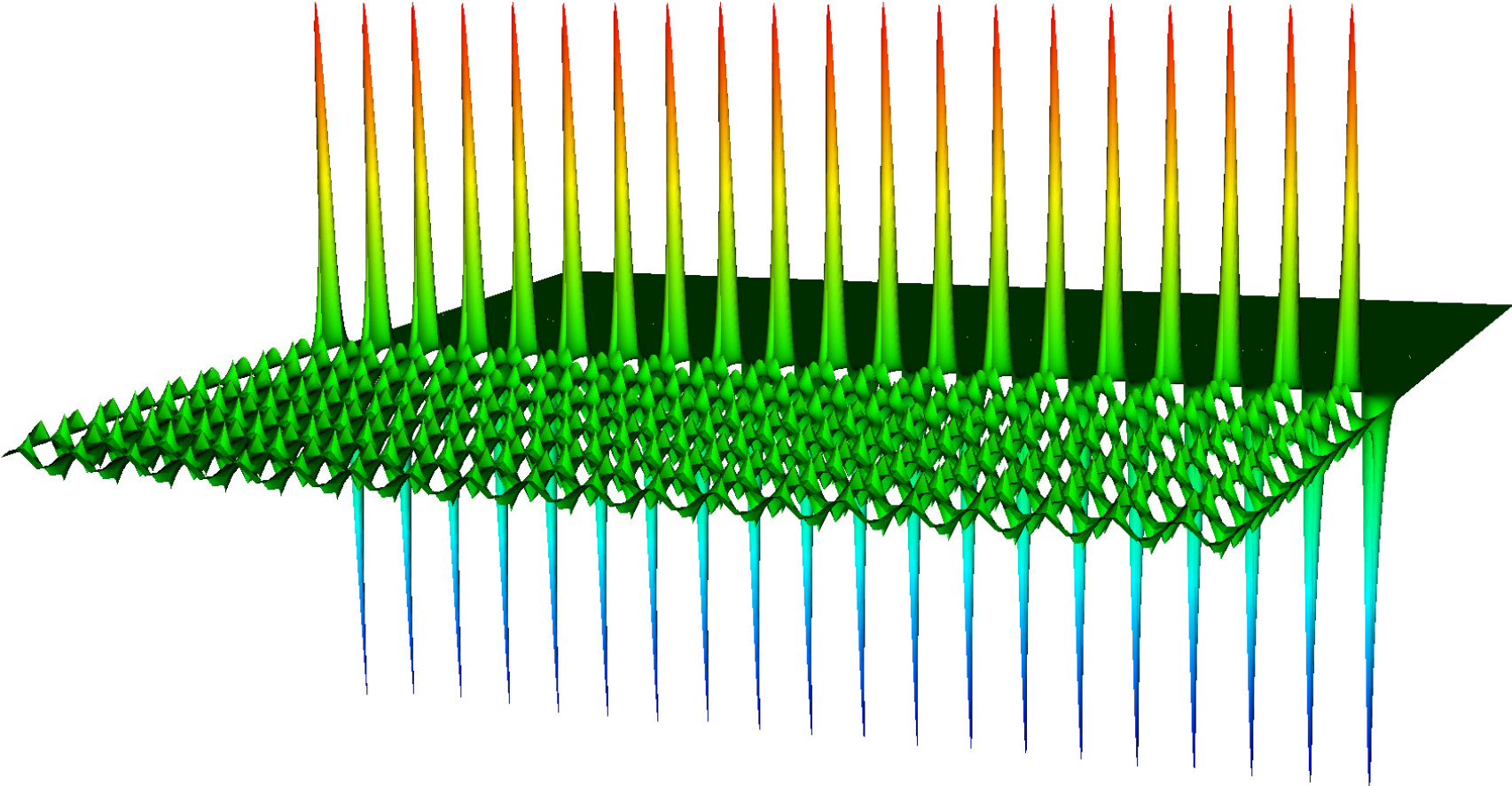}
\caption{Circles as inclusions.}
\end{subfigure}
\begin{subfigure}[b]{0.48\textwidth}
\centering
\includegraphics[width=0.98\textwidth]{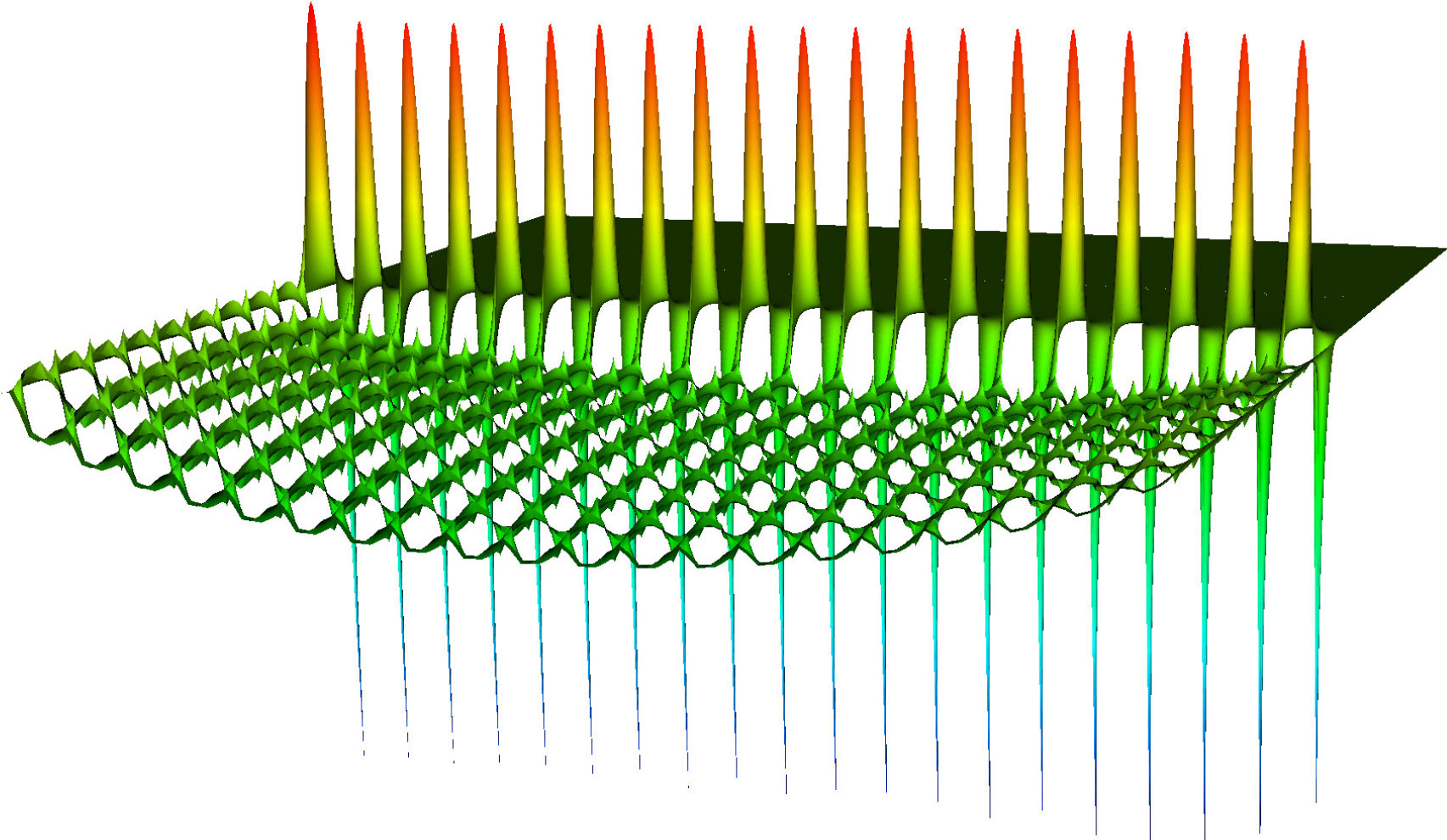}
 \caption{Ellipses as inclusions.}
\end{subfigure}
\caption{Visualization of the pressure with $\epsilon = 10^{-1}$.}
\label{fig.pressure ellipse eps 10^-1}
\end{figure}
These non-periodic boundary conditions introduce a pollution effect in
comparison with the theoretical results valid for the periodic case.
However, we see in figure~\ref{fig.pressure ellipse eps 10^-1} that away
from the boundary a jump in the pressure is visible. For this reason the domain length is set to $L=2$ to reduce the effect of the boundaries.

In analogy to the previous case we define cell-wise average quantities to give the quantitative results shown in table~\ref{tab.eps_con_mean_real}.
The first column of the table shows that the pressure interface condition \eqref{Presspm2A} is not fulfilled by the cell-wise average quantities, while the second column shows that the Beavers-Joseph condition \eqref{4.95} is satisfied.
We observe thus that the pollution effect on the interface condition mainly
concerns the pressure and outer boundary layer effects appear, see \cite{JaegerMN:2001}.
To get rid of this effect we consider the averaged quantities only over part of the interface away from the boundary, i.e. the integral is taken over $\widehat \Gamma := (0.2,1.8)\times \Set 0$ instead of $\Gamma=(0,2)\times \Set 0$.
In this case, as observed in the third and fourth columns, both interface laws are fulfilled, i.e. a convergence with order $\epsilon$ is shown.

\begin{table}
\centering
\begin{tabular}{r|cccc}
\toprule
$\epsilon$&$\norm{\clos{p^\ep_{\operatorname{d}}}-\clos{p^\ep}-\cblw\clos{\partial_2u^\ep_1}}_{2,\Gamma}$&$ \norm{\frac{\clos{u_1^\ep}}{\ep}+\cbl\clos{\partial_2u^\ep_1}}_{2,\Gamma}$&$\norm{\clos{p^\ep_{\operatorname{d}}}-\clos{p^\ep}-\cblw\clos{\partial_2u^\ep_1}}_{2,\widehat\Gamma}$&$ \norm{\frac{\clos{u_1^\ep}}{\ep}+\cbl\clos{\partial_2u^\ep_1}}_{2,\widehat \Gamma}$\\
\midrule
&\multicolumn{4}{c}{\textbf{circular inclusions}}\\
1              &  - &  6.75E-02 &  - &  6.04E-02   \\
$\frac 1 3 $   &  5.11E-11 &  2.25E-02 &  5.03E-11 &  2.01E-02   \\
0.1            &  5.53E-09 &  6.74E-03 &  4.45E-09 &  6.03E-03   \\
$\frac 1 {31}$ &  6.56E-09 &  2.17E-03 &  4.68E-09 &  1.94E-03   \\
0.01           &  1.12E-07 &  6.63E-04 &  6.62E-11 &  5.93E-04   \\
\midrule   
	&\multicolumn{4}{c}{\textbf{oval inclusions}}\\
1                & -  &   5.65E-02  &    -   &   5.06E-02  \\
$\frac 1 3$      & 4.91E-02  &   1.93E-02  &    4.33E-02   &   1.68E-02  \\
0.1              & 4.90E-02  &   6.07E-03  &    2.31E-02   &   4.91E-03  \\
$\frac 1 {31}$   & 3.43E-02  &   2.13E-03  &    9.84E-03   &   1.51E-03  \\
 0.01            & 2.01E-02  &   8.70E-04  &    2.97E-03   &   4.70E-04  \\
\bottomrule 
\end{tabular}
\centering \caption{Numerical evidence of the convergence results on the interface $\Gamma$ in case of non-periodic boundary conditions. }\label{tab.eps_con_mean_real}
\end{table}

\section{Conclusions}
\label{conclusions}
A pressure jump condition of the slow viscous flow over a porous bed has
been rigorously derived by Marciniak-Czochra and Mikeli\'c in a recent
article.
In this work, we have presented a numerical confirmation of this condition
based on finite elements.
A goal oriented mesh adaptivity based on an a-posteriori error estimator has been used to precisely
calculate the needed problems: the Navier boundary layer and some appropriate
cell problems to calculate the permeability tensor.
Two test cases have been shown: a periodic flow and a flow with an
injection condition. The first case is fully supported by the theoretical
results of Marciniak-Czochra and Mikeli\'c and the comparison between upscaled solution and the solution on the microscopic
level confirms the pressure jump law.
The second case introduces a perturbation due to a boundary layer at the inflow and outflow. Numerical
calculations also in this case show the pressure jump away from the
boundaries, where the pollution effect is strongly diminished.
In case of isotropic porous medium, the pressure is continuous due to the
exponential vanishing of the boundary layer pressure. This is confirmed
numerically as well.

\begin{acknowledgments}
AM-C  was supported by ERC Starting Grant "Biostruct" No. 210680 and Emmy Noether Programme of German Research Council (DFG). The research of A.M. was partially supported by the  Programme Inter Carnot Fraunhofer from BMBF (Grant 01SF0804) and ANR. Research visits of A.M. to the Heidelberg University were supported in part by the Romberg professorship at IWR, Heidelberg University, 2011-1013. TC was supported by the German Research Council (DFG) through project
``Modellierung, Simulation und Optimierung der Mikrostruktur mischleitender
SOFC-Kathoden'' (RA 306/17-2).
\end{acknowledgments}

\bibliographystyle{jfm}
\bibliography{refs}
\end{document}